%% file: 1997-027.tex
\documentclass[12pt,draft]{amsart}
\usepackage{amssymb}
\input table.tex %Wichura's macros

\newskip\alignskip \alignskip=10pt

\DeclareMathOperator{\Id}{Id}
\DeclareMathOperator{\Tr}{Tr}
\DeclareMathOperator{\diag}{diag}
\newcommand{\rfrac}[2]{\leavevmode
  \raise.5ex\hbox{$\scriptscriptstyle#1$}\kern-.1em
  {\scriptstyle/}\kern-.1em\lower.2ex\hbox{$\scriptscriptstyle#2$}}

\newcommand{\C}{\mathbb{C}}

\newcommand{\Q}{\mathbb{Q}}

\newcommand{\Z}{\mathbb{Z}}

\renewcommand{\(}{\langle\,}
\renewcommand{\)}{\,\rangle}

\numberwithin{equation}{section}

\newtheorem{prop}[equation]{Proposition}

\newtheorem{lemma}[equation]{Lemma}

\theoremstyle{definition}

\newcounter{parts}
\def\prt(#1){{\rm({\it#1\/})}}
\newenvironment{itlist}%
  {\begin{list}%
    {\hbox to 6pt{\hfil}{\rm(\alph{parts})}}%
    {\usecounter{parts}\setlength{\labelwidth}{24pt}%
     }%
  }%
  {\end{list}}

\begin{document}
\title[Representations of the Iwahori-Hecke Algebra $HF_4$]{Explicit
Irreducible Representations of the Iwahori-Hecke Algebra of Type $F_4$}
\author{Arun Ram}
\address{\hskip-\parindent
Arun Ram\\Department of Mathematics\\Princeton University\\
Princeton, NJ 08544-1000}
\email{rama@math.princeton.edu}
\author{D. E. Taylor}
\address{\hskip-\parindent
D. E. Taylor\\School of Mathematics and Statistics\\
University of Sydney\\ NSW 2006 Australia}
\email{don@maths.su.oz.au}
\thanks{
The research of the first author was
supported in part by an Australian Research Council
Research Fellowship, National Science Foundation grant DMS-9622985,
and the Mathematical Sciences Research Institute.
Research at MSRI is supported in part by NSF grant DMS-9022140.}
%\date{24 March 1997}
\begin{abstract}
A general method for computing irreducible representations of Weyl
groups and Iwahori-Hecke algebras was introduced by the first author in
\cite{ram:1997}. In that paper the representations of the algebras of
types $A_n$, $B_n$, $D_n$ and $G_2$ were computed and it is the purpose
of this paper to extend these computations to $F_4$. The main goal here
is to compute irreducible representations of the Iwahori-Hecke algebra
of type $F_4$ by only using information in the character table of the
Weyl group.
\end{abstract}
\maketitle

\section{Introduction}

In his thesis \cite{hoefsmit:1974} P.N. Hoefsmit wrote down explicit
irreducible representations of the Iwahori-Hecke algebras 
$HA_{n-1}$, $HB_n$, and $HD_n$, of
types $A_{n-1}$, $B_n$, and $D_n$, respectively.  
Hoefsmit's thesis was never
published and H. Wenzl \cite{wenzl:1988}
independently discovered these representations
in the type $A_{n-1}$ case.
The irreducible representations of Hoefsmit are analogues of the
``seminormal'' representations of the Weyl groups of types $A_{n-1}$,
$B_n$ and $D_n$ which were written down by A. Young \cite{young:1932}.
The Iwahori-Hecke algebras depend on parameters $p$ and $q$ and
one can recover the representations of Young by setting
$p$ and $q$ equal to 1 in Hoefsmit's representations.

In this paper we shall extend Hoefsmit's result 
and determine explicit realizations of all the irreducible
representations of the Iwahori-Hecke
algebra $HF_4$. The matrix entries of these representations
are well defined when $p=q=1$
and, when one sets $p=q=1$, our representations
specialize to give explicit realizations of all the irreducible
representations of the Weyl group of type $F_4$.
The final results are tabulated in the last section of this paper.
Our numbering scheme for the
irreducible characters of $HF_4$ follows Geck~\cite{geck:1994}.

In analogy with Hoefsmit, our representations of $HF_4$ are in 
``seminormal'' form with respect to the chain of subalgebras
\[
HF_4\supseteq HB_3\supseteq HA_2\supseteq HA_1,
\] 
which means that we choose the irreducible representations $\varphi^k$
of $HF_4$ so that, for all $h\in HB_3$, the matrices $\varphi^k(h)$ 
are block diagonal matrices where the blocks $\varphi^{\mu}(h)$
are determined by the irreducible representations $\varphi^\mu$
of $HB_3$.
Similarly, for all $h\in HA_2$, the matrices $\varphi^k(h)$ 
are block diagonal where the
blocks are determined by the irreducible representations $\varphi^\lambda(h)$
of $HA_2$.
In this way we construct the irreducible representations of
$HF_4$ inductively, by using the
branching rules for restricting
representations from $HF_4$ to $HB_3$ and from $HB_3$ to $HA_2$.
These branching rules can be calculated easily from the character
tables of the corresponding Weyl groups.

The matrix entries of our representations are
rational functions in the variables $p$ and $q$.
These rational functions are quotients of polynomials
in $\Z[p,q,p^{-1},q^{-1}]$ and the denominators
contain only the polynomials
\begin{equation}\label{eqn:badvalues}
[2]_p,\enspace [2]_q,\enspace 
[2]_{pq},\enspace [2]_{pq^{-1}},\enspace 
[2]_{p^2q},\enspace [2]_{p^2q^{-1}},\enspace
[3]_p,\quad\text{and}\quad [3]_q,
\end{equation}
where $[2]_x = x+x^{-1}$ and $[3]_x=x^2+1+x^{-2}$.
This means that our representations are well defined over any
field $F$ such that $p,q\in F$ and none of the polynomials
in (\ref{eqn:badvalues}) are equal to $0$.

The research in this paper was begun 
while the first author was visiting Sydney University on 
an Australian Research Council Research fellowship.  
Arun Ram thanks Sydney University and
especially G. Lehrer for their wonderful hospitality
during the year that he spent in Sydney.  The 
first author is also grateful for support 
from a postdoctoral fellowship at the
Mathematical Sciences Research Institute in Berkeley
where the writing of this paper was completed.

\section{Preliminaries}

Let $p$ and $q$ be indeterminates.
The Iwahori-Hecke algebra $HF_4$ 
is the associative algebra with 1 over
the field $\C(p,q)$ generated by $T_1, T_2, T_3, T_4$ with relations
\begin{align*}
T_1T_2T_1&=T_2T_1T_2\\
T_3T_4T_3&=T_4T_3T_4\\
T_2T_3T_2T_3&=T_3T_2T_3T_2\\
T_iT_j&=T_jT_i,\qquad\text{if $j\ne i\pm 1$},\\
T_i^2&=(p-p^{-1})T_i+1,\quad\text{for $i=1,2$,}\\
T_i^2&=(q-q^{-1})T_i+1,\quad\text{for $i=3,4$}.
\end{align*}
This is the Iwahori-Hecke algebra corresponding to the Weyl group $WF_4$.
The Weyl group $WF_4$ is generated by $s_1,s_2,s_3,s_4$ which satisfy the same
relations as the $T_i$ except with $p=q=1$.  Let $HA_1$, $HA_2$, and $HB_3$
be the subalgebras of $HF_4$ such that
\begin{align*}
HA_1 & 
\quad\text{is generated by $T_1$,} \\
HA_2 & 
\quad\text{is generated by $T_1$ and $T_2$, and} \\
HB_3 & 
\quad\text{is generated by $T_1$, $T_2$ and $T_3$.} 
\end{align*}
These are the Iwahori-Hecke algebras corresponding to the Weyl groups
$WA_1=<s_1>$, $WA_2 = <s_1,s_2>$ and $WB_3=< s_1,s_2,s_3>$,
respectively.

Our goal in this paper is to compute explicit representations of 
$HF_4$
using only the information in the character tables of the Weyl groups
$WA_1$, $WA_2$, $WB_3$ and $WF_4$.
We shall use the following notations.
\begin{itlist}
\item\notag $d_\lambda$ will denote the dimension of the irreducible
representation  indexed by $\lambda$;
\item\notag $\chi^\lambda$ will denote the character of the
irreducible representation of the Weyl group $W$ indexed by $\lambda$;
\item\notag $\Id_\lambda$ will denote the $d_\lambda\times d_\lambda$
identity matrix;
\item\notag $T_w$, $w\in W$, will denote the usual basis of the
Iwahori-Hecke algebra $H$ given by
$T_w=T_{i_1}\cdots T_{i_k}$ if $w=s_{i_1}\cdots s_{i_k}$ is a reduced word
for $w$.
\item\notag  If $A$ and $B$ are matrices then
$A\oplus B$ and $A\otimes B$
will denote the standard operations of direct sum and
tensor product of matrices.
\end{itlist}

\noindent
We shall need the following well known facts:

\medskip\noindent
FACT 1.  The irreducible representations of the Iwahori-Hecke 
algebra are indexed in the
same way as the corresponding Weyl group.  Thus,
\begin{itlist}
\item
The irreducible representations of $HF_4$ are indexed by 
$k\in \{1,2,\ldots,25\}$ (in the same manner as in \cite{geck:1994}
and in the same order as in the table on p. 412 of \cite{carter:1985}).
\item
The irreducible representations of $HB_3$ are indexed by 
pairs of partitions $(\alpha,\beta)$ such that $|\alpha|+|\beta| = 3$.
\item
The irreducible representations of $HA_2$ are indexed by 
partitions $\lambda$ of $3$.
\item
The irreducible representations of $HA_1$ are indexed by 
partitions $\gamma$ of $2$.
\end{itlist}

\medskip\noindent
FACT 2. \cite{carter:1985} \S 10.11.
The dimension of an irreducible
Iwahori-Hecke algebra representation is the same as that of the corresponding
representation of the Weyl group and the branching rules for Iwahori-Hecke 
algebras are the same as for the corresponding Weyl groups.
Thus the branching rules for the inclusions $HF_4\supseteq HB_3\supseteq HA_2$
can be calculated directly from the character tables of the 
corresponding Weyl groups.
We have tabulated these branching rules
in Tables~\ref{tbl:1} and \ref{tbl:2}.  

\medskip\noindent
FACT 3. \cite{carter:1985} \S 10.9 and \cite{curtis:1981} (9.21).  
Let $H$ be an Iwahori-Hecke algebra and let $W$ be the 
corresponding Weyl group.  
If $\lambda$ is an index for an irreducible representation of the Iwahori-Hecke
algebra $H$ then the minimal central idempotent corresponding to $\lambda$
can be written in the form
\[
z_\lambda = \sum_{w\in W} z^\lambda_w T_w,
\]
where $z^\lambda_w\in \C(p,q)$ are elements which are well defined when
$p=q=1$.  Furthermore, at $p=q=1$,
\begin{equation}\label{eqn:idem}
z_\lambda\big|_{p=q=1} =
 \frac{\chi^\lambda(1)}{|W|}\sum_{w\in W} \chi^\lambda(w^{-1})w,
\end{equation}
where $\chi^\lambda$ is the character of the irreducible
representation of $W$ indexed by $\lambda$.

\medskip\noindent
FACT 4.  
Let $H$ be an Iwahori-Hecke algebra and let $W$ be the 
corresponding Weyl group.  
Let $R$ be the root system corresponding to $W$ and
let 
\begin{align*}
r_s &= \text{a reflection in a short root,} \\
r_l &= \text{ a reflection in a long root,}\\
N_s &= \text{ the number of positive short roots in $R$, and}\\
N_l &= \text{ the number of positive long roots in $R$.} 
\end{align*}
If there is only one root length then we declare all roots to be
short.
For each $\lambda$ indexing an irreducible representation of $H$
let $\chi^\lambda$ be the character of the corresponding
irreducible representation of the Weyl group and define
\begin{equation}\label{eqn:central}
c(\lambda) =  \chi^\lambda(w_0)
p^{c(\lambda,s)}q^{c(\lambda,\ell)} 
\end{equation}
where
\[
c(\lambda,s) = \frac{N_s\chi^\lambda(r_s)}{\chi^\lambda(1)}
\quad\text{and}\quad
c(\lambda,l) = \frac{N_l\chi^\lambda(r_l)}{\chi^\lambda(1)}.
\]
Let $\varphi^\lambda$ be a realization of the irreducible representation
indexed by  $\lambda$ and let $\Id_\lambda$ be the 
$d_\lambda\times d_\lambda$ identity matrix, where $d_\lambda$ is the
dimension of $\varphi^\lambda$.  Then we have the following result
\cite{kilmoyer:1969}, \cite{geck:1994}, \cite{ram:1997}:
\begin{itlist}
\item If $w_0$ is central in $W$ then 
$\varphi^\lambda(T_{w_0}) = c(\lambda)\Id_\lambda,$ 
\item If $w_0$ is not central in $W$ then
$\varphi^\lambda(T_{w_0}^2) = c(\lambda)^2\Id_\lambda.$ 
\end{itlist}

\section{Seminormal representations}

We shall compute the irreducible representations of $HF_4$ 
inductively: the representations of $HA_1$ are
one dimensional and one can immediately write them down,
then we compute irreducible representations
of $HA_2$, then $HB_3$, and finally $HF_4$.
At each step we use the information from the previous 
cases since we construct the representations
such that upon restriction to any of these subalgebras
they are in block diagonal form with diagonal blocks determined 
by the previous calculations.
The irreducible representations of $HA_2$ are easy to derive
and the irreducible representations of 
$HB_3$ can be derived in a similar fashion to the
way that we complete the calculations for $HF_4$ below.
Thus, in our description below we shall {\it assume that
the irreducible representations of $HA_1$, $HA_2$, and
$HB_3$ are already known} and we shall describe how to obtain
the irreducible representations of $HF_4$.
The irreducible ``seminormal'' representations of $HA_1$,
$HA_2$, and $HB_3$ are tabulated in 
Section~\ref{secn:branch}
below.

Let $k$ be an index for an irreducible representation of $HF_4$.
The branching rule 
$$\varphi^{k}\downarrow_{HB_3} \cong
\varphi^{\mu^{(1)}} \oplus \varphi^{\mu^{(2)}}\oplus \cdots 
\oplus \varphi^{\mu^{(\ell)}}$$
describing the restriction of representations of
$HF_4$ to $HB_3$ can be computed from the character table of
the corresponding Weyl groups.
We shall say that the irreducible representation $\varphi^k$ of $HF_4$ is
in {\it seminormal form} if 
\begin{equation}\label{eqn:sncondition}
\varphi^{k}(h) =
\varphi^{\mu^{(1)}}(h) \oplus \varphi^{\mu^{(2)}}(h) \oplus \cdots 
\oplus \varphi^{\mu^{(\ell)}}(h),
\qquad\text{for all $h\in HB_3$.}
\end{equation}
We require the two sides of (\ref{eqn:sncondition}) to be equal
as matrices.

We shall compute irreducible representations of $HF_4$
which are in seminormal form.
Assuming that the irreducible representations of
$HB_3$ are known, the seminormal condition 
implies that to determine the irreducible
representations of $HF_4$ it is only necessary to determine 
the matrices $\varphi^k(T_4)$ for each $k$.

Suppose that $\varphi^k$ and $\psi^k$ are two solutions to this problem, 
i.e.  $\varphi^k$ and
$\psi^k$ are both realizations of the irreducible representation of 
$HF_4$ indexed
by $k$ and we have
$$\varphi^{k}(h) =
\psi^{k}(h) =
\varphi^{\mu^{(1)}}(h) \oplus \varphi^{\mu^{(2)}}(h) \oplus \cdots 
\oplus \varphi^{\mu^{(\ell)}}(h),
$$
for all $h\in HB_3$.
Then there is a matrix $P\in GL(d_k)$, where $d_k$ is the 
dimension of $\varphi^k$, such that 
$P\varphi^{k}(h)P^{-1} = \psi^{k}(h)$,
for all $h\in HF_4$.  By Schur's lemma this matrix is unique up to constant
multiples.  On the other hand we have
\begin{align*}
P
(\varphi^{\mu^{(1)}}(h) \oplus \varphi^{\mu^{(2)}}(h) \oplus \cdots 
\oplus \varphi^{\mu^{(\ell)}}(h)) P^{-1} 
=P\varphi^{k}(h)P^{-1} 
= \psi^{k}&(h) \\
=\varphi^{\mu^{(1)}}(h) \oplus \varphi^{\mu^{(2)}}(h) \oplus \cdots 
\oplus \varphi^{\mu^{(\ell)}}&(h), 
\end{align*}
for all $h\in HB_3$.  By inspection of the 
table of branching rules from $HF_4$ to $HB_3$
one sees that the summands  $\varphi^{\mu^{(i)}}$
are all distinct irreducible representations of $HB_3$. 
Hence, Schur's lemma implies that
\begin{equation}\label{eqn:trans}
P = p_1\Id_{\mu^{(1)}}\oplus p_2\Id_{\mu^{(2)}}\oplus 
\cdots \oplus p_\ell\Id_{\mu^{(\ell)}},
\end{equation}
where the $p_i$ are nonzero constants.  Replacing $P$ by $p_1^{-1}P$
we may suppose that $p_1=1$.
Conversely, any choice of 
$p_i\ne 0$, $p_1=1$, in the equation (\ref{eqn:trans})
defines a matrix $P$ such that 
$P\varphi^{k}P^{-1}$ is a seminormal
representation.  Thus we have the following result.

\begin{prop}\label{prop:param}
If $\varphi^k$ is in seminormal form then
the matrix $\varphi^{k}(T_4)$ is determined up to the choice of
$\ell-1$ free parameters where $\ell$ is the number of irreducible summands
in $\varphi^k$ on restriction to $HB_3$.
\end{prop}

Let $w_{0,1}$, $w_{0,2}$, $w_{0,3}$, and $w_{0,4}$
be the longest elements in the Weyl groups $WA_1$, $WA_2$,
$WB_3$, and $WF_4$, respectively.
Define elements
\begin{align}
D_1 &= T_{w_{0,1}} = T_1, \notag\\
D_2 &= T_{w_{0,2}}^2 = (T_1T_2T_1)^2, \notag\\
D_3 &= T_{w_{0,3}} = (T_3T_2T_1)^3, \label{eqn:two}\\
D_4 &= T_{w_{0,4}} = (T_4T_{w_{0,3}})^3 T_{w_{0,2}}^{-2}.\notag 
\end{align}
in $HF_4$.

\begin{lemma}\label{lemma:djdet}
If $\varphi^k$ is in seminormal form then
the matrices $\varphi^k(D_j)$ are uniquely determined,
for all $1\le k\le 25$, $1\le j\le 4$.
\end{lemma}
\begin{proof}
This follows from Fact 4.
\end{proof}

\noindent
The matrices $\varphi^k(D_j)$ are tabulated in \ref{secn:f4rep}.

Let $\sigma$ be a permutation matrix such that
\begin{align}
\sigma \varphi^k(h) \sigma^{-1}
&= \varphi^{(1^3)}(h)\oplus\cdots \oplus \varphi^{(1^3)}(h) \notag\\
&\phantom{= \varphi^{(1^3)}(h)\oplus\cdots }
\oplus\enspace
\varphi^{(21)}(h)\oplus\cdots \oplus \varphi^{(21)}(h) \notag\\
&\phantom{= \varphi^{(1^3)}(h)\oplus\cdots 
\oplus\enspace
\varphi^{(21)}(h)}
\oplus\enspace
\varphi^{(3)}(h)\oplus\cdots \oplus \varphi^{(3)}(h) \label{eqn:a2mats}\\
&=\bigoplus_{\lambda\vdash 3} \varphi^\lambda(h)^{\oplus m_\lambda} \notag
\end{align}
for all $h\in HA_2$.  
The constant $m_\lambda$ is the number of times
the matrix $\varphi^\lambda(h)$ appears.
Since $\sigma \varphi^k(T_4)\sigma^{-1}$ commutes with all of the
matrices in (\ref{eqn:a2mats}), it follows from Schur's lemma that 
\begin{align*}
\sigma \varphi^k(T_4) \sigma^{-1}
&= (T^k_{(3)}\otimes \Id_{(3)}) 
\oplus (T^k_{(21)}\otimes \Id_{(21)})
\oplus (T^k_{(1^3)}\otimes\Id_{(1^3)}) \\
&= \bigoplus_{\lambda\vdash 3} T^k_\lambda\otimes \Id_\lambda, 
\end{align*}
where, for each $\lambda$, $T^k_\lambda$ is an 
$m_\lambda\times m_\lambda$ matrix
and $\Id_\lambda$ is the $d_\lambda\times d_\lambda$ identity matrix.
Note that 
\begin{equation}\label{eqn:tkmin}
T^k_\lambda\otimes \Id_\lambda 
=\sigma \varphi^k(z_\lambda T_4)\sigma^{-1},
\end{equation}
where $z_\lambda$ is the minimal central idempotent in $HA_2$
corresponding to $\lambda$.
We can use the same method to write
\[
\sigma\varphi^k(D_3)\sigma^{-1}
=\bigoplus_{\lambda\vdash 3} D^k_\lambda\otimes \Id_\lambda,
\]
where $D_3=T_{w_{0,3}}$, as given in (\ref{eqn:two}).

To determine the matrices $\varphi^k(T_4)$ it is sufficient to
determine the matrices $T^k_\lambda$.  The matrices
$D^k_\lambda$ are completely determined by Lemma \ref{lemma:djdet}
and can easily be determined from the tables in \ref{secn:f4rep}.
The relations $(T_4D_3)^3=D_4D_2^2$  
and the relation $T_4^2 = (q-q^{-1})T_4+1$ imply that
\begin{equation}\label{eqn:relations}
(T^k_\lambda D^k_\lambda)^3 = 
c(k)c(\lambda)^2\Id_{m_\lambda}
\quad\text{and}\quad
(T^k_\lambda)^2 = (q-q^{-1})T^k_\lambda +\Id_{m_\lambda},
\end{equation}
where $c(k)$ and $c(\lambda)$ are the constants given in 
equation (\ref{eqn:central}).

\subsection{Determining the diagonal entries of $\varphi^k(T_4)$}

We shall determine the diagonal entries of the matrices the matrices
$T^k_\lambda$ by determining the traces of the matrices
\[
T^k_\lambda(D^k_\lambda)^{-2}, \quad
T^k_\lambda (D^k_\lambda)^{-1}, \quad
T^k_\lambda,  \quad 
T^k_\lambda D^k_\lambda,
\quad\text{and}\quad
T^k_\lambda (D^k_\lambda)^2.
\]

\begin{prop}\label{prop:trace}  
Fix an index $k$ for an irreducible representation of
$HF_4$ and let $\lambda$ be an index for an irreducible
representation of $HA_2$. Let $T^k_\lambda$, $D^k_\lambda$ 
and $z_\lambda$ be as above and
let $\chi^k$ and $\chi^\lambda$ be the irreducible characters of the
Weyl groups $WF_4$ and $WA_2$ which correspond to $k$ and $\lambda$,
respectively. Let $c(k)$ and $c(\lambda)$ be
the constants defined in (\ref{eqn:central}).  Then
\begin{align*}
\Tr(T^k_\lambda) &= \frac{1}{12}\sum_{w\in WA_2} \chi^\lambda(w^{-1})
\left((q-q^{-1})\chi^k(w)+(q+q^{-1})\chi^k(ws_4)\right),
\tag{a}\\
\Tr(T^k_\lambda D^k_\lambda) 
&= \frac 
{\chi^k(w_{0,4})c(k)^{\frac{1}{3}}c(\lambda)^{\frac{2}{3}} }
{6}
\sum_{w\in WA_2} \chi^\lambda(w^{-1})\chi^k(ws_4w_{0,3}),
\tag{b}\\
\Tr((T^k_\lambda D^k_\lambda)^2) &= 
\frac 
{c(k)^{\frac{2}{3}}c(\lambda)^{\frac{4}{3}} }
{6}
\sum_{w\in WA_2}\chi^\lambda(w^{-1})\chi^k(w(s_4w_{0,3})^2).
\tag{c}
\end{align*}
\end{prop}
\begin{proof}
(a)
From the second equation in (\ref{eqn:relations}) we have that
each eigenvalue of 
$T^k_\lambda$ is either $q$ or
$-q^{-1}$ and consequently $\Tr(T^k_\lambda)=t_1q-t_2q^{-1}$ for 
some positive integers $t_1$ and $t_2$.  
These constants are determined as follows.
Using (\ref{eqn:tkmin}) we get that
\begin{equation}\label{eqn:parta}
\begin{align*}
t_1-t_2 &=
\Tr(T^k_\lambda)\big|_{p=q=1} 
= \frac{1}{\chi^\lambda(1)}\Tr(T^k_\lambda\otimes \Id_\lambda)
\big|_{p=q=1}\\
&= \frac{1}{\chi^\lambda(1)}\Tr(\sigma\varphi^k(z_\lambda T_4)\sigma^{-1})
\big|_{p=q=1}
=\frac{1}{\chi^\lambda(1)}\Tr(\varphi^k(z_\lambda T_4))\big|_{p=q=1}. 
\end{align*}
\end{equation}
Then we use (\ref{eqn:idem}) to obtain
%\begin{align*}
\[
t_1-t_2 =\frac{1}{\chi^\lambda(1)}
\chi^k\left(\frac{\chi^\lambda(1)}{6}
\sum_{w\in WA_2} \chi^\lambda(w^{-1})ws_4\right) 
=\frac{1}{6}\sum_{w\in WA_2} \chi^\lambda(w^{-1})\chi^k(ws_4). 
\]
%\end{align*}
If $\Id^k_\lambda$ is the identity matrix of the 
same dimension as $T^k_\lambda$ then
\begin{align*}
t_1+t_2 
&= \Tr(\Id^k_\lambda) \big|_{p=q=1}
= \frac{1}{\chi^\lambda(1)} 
\Tr(\Id^k_\lambda\otimes \Id_\lambda) \big|_{p=q=1} \\
&= \frac{1}{\chi^\lambda(1)} \Tr(\varphi^k(z_\lambda))\big|_{p=q=1} 
=\frac{1}{6}\sum_{w\in WA_2} \chi^\lambda(w^{-1})\chi^k(w).
\end{align*}
These two equations determine $t_1$ and $t_2$ and thus $\Tr(T^k_\lambda)$
is determined.

(b)
It follows from Fact 4 and the first equation in
(\ref{eqn:relations}) that the eigenvalues of $T^k_\lambda D^k_\lambda$ 
are all of the form
$\omega^i c(\lambda)^{\frac{2}{3}} c(k)^{\frac{1}{3}}$
where $\omega$ is a primitive cube root of unity.  Hence
\[
\Tr(T^k_\lambda D^k_\lambda) 
= \eta  c(\lambda)^{\frac{2}{3}} c(k)^{\frac{1}{3}}
\]
for some constant $\eta\in \C$.  By setting $p$ and $q$ equal to $1$
we have
\begin{align*}
\eta \chi^k(w_{0,4})
&=
\Tr(T^k_\lambda D^k_\lambda)\big|_{p=q=1} 
= \frac{1}{\chi^\lambda(1)}
\Tr(T^k_\lambda D^k_\lambda\otimes \Id_\lambda)\big|_{p=q=1} \\
&=  \frac{1}{\chi^\lambda(1)}
\Tr(\varphi^k(z_\lambda T_4T_{w_{0,3}}))\big|_{p=q=1}
\end{align*}
as in (\ref{eqn:parta}).  Using (\ref{eqn:idem}) we get
\[
\eta\chi^k(w_{0,4})
= \frac{1}{6}
\sum_{w\in WA_2} \chi^\lambda(w^{-1})\chi^k(ws_4w_{0,3}).
\]

The proof of (c) is similar to that of (b) once one notes that
Fact 4 and the first equation in (\ref{eqn:relations}) imply that
the eigenvalues of the matrix $(T^k_\lambda D^k_\lambda)^2$ 
are all of the form
$\omega^{2i} c(\lambda)^{\frac{4}{3}} c(k)^{\frac{2}{3}}$.
\end{proof}

\begin{lemma}\label{lemma:trace}
Given matrices $T$ and $D$ such that
$T^2=(q-q^{-1})T+\Id$ and $(TD)^3 = c\Id$ where $c$ is a constant, we
have
\begin{itlist}
\item
$\Tr(TD^{-1}) = (q-q^{-1})\Tr(D^{-1}) + c^{-1}\Tr( (TD)^2)$.
\item
$\Tr(TD^2) = c\Tr(D^{-1}) - (q-q^{-1}) \Tr( (TD)^2)$.
\item
$\Tr(TD^{-2}) = (q-q^{-1})\Tr(D^{-2}) + c^{-1}(q-q^{-1})\Tr(TD)
+c^{-1}\Tr(D)$.
\end{itlist}
\end{lemma}

\begin{proof}
(a)\enspace Writing the given equations in the form $T = (q-q^{-1})\Id +
T^{-1}$ and $(TD)^{-1} = c^{-1} (TD)^2$, we have
\begin{align*}
\Tr(TD^{-1})
&= (q-q^{-1})\Tr(D^{-1}) + \Tr(T^{-1}D^{-1})\\
&= (q-q^{-1})\Tr(D^{-1}) + c^{-1}\Tr(TDTD).
\end{align*}
(b)\enspace Similarly, from the fact that $T^{-2} = \Id -
(q-q^{-1})T^{-1}$,
\begin{align*}
\Tr(TD^{2})
&= \Tr(DTD) = c\Tr(T^{-1}D^{-1}T^{-1}) = c\Tr(T^{-2}D^{-1})\\
&=c\Tr(D^{-1})-c(q-q^{-1})\Tr(T^{-1}D^{-1}) \\
&=c\Tr(D^{-1})-(q-q^{-1})\Tr(TDTD). 
\end{align*}
\begin{align}
\Tr(TD^{-2})
&=(q-q^{-1}) \Tr(D^{-2}) + \Tr(T^{-1}D^{-2})\tag{c}\\
&=(q-q^{-1}) \Tr(D^{-2}) + \Tr(D^{-1}T^{-1}D^{-1})\notag \\
&=(q-q^{-1}) \Tr(D^{-2}) + c^{-1} \Tr(TDT) \notag\\
&=(q-q^{-1}) \Tr(D^{-2}) + c^{-1} \Tr(T^2D) \notag\\
&=(q-q^{-1}) \Tr(D^{-2}) + c^{-1}(q-q^{-1}) \Tr(TD) + c^{-1}\Tr(D). \notag
\end{align}
\end{proof}

Assume that $T^k_\lambda$ has dimension at most $5$
and write $D^k_\lambda = \diag(d_1,d_2,\ldots,d_r)$.
The diagonal entries of $D^k_\lambda$ are determined by Proposition
\ref{lemma:djdet} and one can check directly that these 
diagonal entries are always all distinct.  
Let $S$ be a subset of 
$\{1,2,\ldots,r\}\backslash \{i\}$ such that $S$ and its complement have at most
$2$ elements.  
Then the diagonal entries of $T^k_\lambda$ are given by
\begin{align}
(T^k_\lambda)_{ii} &= \Tr(TE_{ii})
\quad\text{where,\quad for each $1\le i\le r$,}\label{eqn:diags}\\
E_{ii} &=
\biggl(\prod_{\substack{j\in S\\j\ne i}}
	   \dfrac{D^k_\lambda-d_j}{d_i-d_j} \biggr)
\biggl(\prod_{\substack{j\not\in S\\j\ne i}}
	   \dfrac{(D^k_\lambda)^{-1}-d_j^{-1}}{d_i^{-1}-d_j^{-1}} \biggr).
   \notag
\end{align}
These values can be evaluated explicitly by expanding $E_{ii}$ in 
terms of $(D^k_\lambda)^j$
and using Lemma~\ref{lemma:trace} and Proposition~\ref{prop:trace}
to evaluate the traces
$\Tr(T^k_\lambda (D^k_\lambda)^j)$.

Formula (\ref{eqn:diags}) suffices for computing the 
diagonal entries of the matrices 
$T^k_\lambda$, and thus of the matrices
$\varphi^k(T_4)$, for all $k$ except $k=25$.  
The matrix $T^{25}_{(21)}$ has dimension
$6$ and formula (\ref{eqn:diags}) is not applicable.
The diagonal entries of the matrix $\varphi^{25}(T_4)$ are computed as follows.
Since the matrices $T^{25}_{(1^3)}$ and $T^{25}_{(3)}$ are each
of dimension two
we use formula (\ref{eqn:diags}) to determine their diagonal entries.
By Lemma~\ref{lemma:trace} and Proposition~\ref{prop:trace}
we can determine the traces of the matrices
\[
T^{25}_{(21)}(D^{25}_{(21)})^{-2}, \quad
T^{25}_{(21)}(D^{25}_{(21)})^{-1}, \quad
T^{25}_{(21)},\quad  
T^{25}_{(21)}D^{25}_{(21)}, \quad
T^{25}_{(21)}(D^{25}_{(21)})^2,
\]
and these traces give five linear relations that the diagonal 
entries of $T^{25}_{(21)}$
must satisfy.  Finally, we use the formula
$$0=\Tr(\varphi^{25}(T_4T_3T_2T_1)) 
= \sum_i
\varphi^{25}(T_4)_{ii}\varphi^{25}(T_3)_{ii}
\varphi^{25}(T_2)_{ii}\varphi^{25}(T_1)_{ii}$$
to determine the diagonal entries of $\varphi^{25}(T_4)$ completely.
This last formula is a consequence of the following lemma.

\begin{lemma}\label{lemma:coxeter}
\begin{itlist}
\item $\Tr(\varphi^{25}(T_4T_3T_2T_1))=0$.
\item The diagonal entries of the matrix $\varphi^k(T_4T_3T_2T_1)$
satisfy
\[\varphi^{k}(T_4T_3T_2T_1)_{ii} = 
\varphi^{k}(T_4)_{ii}\varphi^{25}(T_3)_{ii}
\varphi^{k}(T_2)_{ii}\varphi^{25}(T_1)_{ii}.
\]
\end{itlist}
\end{lemma}
\begin{proof}
(a)  Since the Coxeter number for the Weyl group $WF_4$ is
$12$ (see \cite{bourbaki:1968}) we have that $(T_4T_3T_2T_1)^6 = T_{w_{0,4}}$. 
Then it follows from Fact 4 that the eigenvalues of the 
matrix $\varphi^{25}(T_4T_3T_2T_1)$ must be of the form
$\omega^j c(25)$ where $\omega$ is a primitive $6$th root of unity
and $c(25)$ is the constant given in (\ref{eqn:central}).
It follows that $\Tr(\varphi^{25}(T_4T_3T_2T_1))
=\eta c(25)$ for some constant
$\eta$.  Then, from the character table of the Weyl group $WF_4$, we have
$$\eta = \eta c(25)\big|_{p=q=1} 
=\Tr(\varphi^{25}(T_4T_3T_2T_1))\big|_{p=q=1}
=\chi^{25}(s_4s_3s_2s_1)
=0.$$
(b)  Let $h\in HB_3$. 
Let $z_\lambda$, $\lambda\vdash 3$, be the minimal
central idempotents in $HA_2$.
Since $\varphi^k$ is in seminormal form the matrices 
$\varphi^k(z_\lambda)$ are diagonal matrices
with $1$'s and $0$'s on the diagonal, their sum is the identity matrix
and they are mutually orthogonal. It follows that there
is a single partition $\lambda$ such that 
\[
\varphi^k(T_4h)_{ii} = 
\varphi^k(z_\lambda T_4h)_{ii}
= \big( \varphi^k(z_\lambda T_4)\varphi^k(z_\lambda h)\big)_{ii}
\]
It follows from the seminormal condition and the fact that
the branching rules for restricting representations from $HF_4$
to $HB_3$ are multiplicity free that the matrix
$\varphi^k(z_\lambda h)$ is a diagonal matrix.  Thus
the diagonal entries of the matrix $\varphi^k(T_4h)$ satisfy
\[
\varphi^k(T_4 h)_{ii} = \varphi^k(T_4)_{ii}\varphi^k(h)_{ii},
\]
for all $h\in HB_3$.
Since the irreducible representations of $HB_3$ and $HA_2$ that we
are using are also chosen to be in seminormal form, their
representations also satisfy a similar identity.  
The result then follows by induction.
\end{proof}

\subsection{Computing the off-diagonal entries of $\varphi^k(T_4)$}

\begin{prop}\label{prop:linear}
Let $T = (t_{ij})$ and  $D = \diag(d_1,d_2,\dots,d_r)$ be $r\times r$
matrices such that the diagonal entries of $D$ are distinct,
none of the entries $t_{ij}$ of $T$ are 0, and 
\[
T^2 = (q-q^{-1})T + \Id\quad\text{and}\quad (TD)^3 = c\Id,
\] 
for some constant $c$.
For distinct indices $i,j,k$ define $u_{ij} = u_{ji} =  t_{ij}t_{ji}$
and $v_{ijk} = t_{ij}t_{jk}/t_{ik}$.  Then
\begin{align}
\sum_{j\ne i} u_{ij} &= -t_{ii}^2 + (q-q^{-1})t_{ii} + 1,
\tag{a}\\
\sum_{j\ne i} u_{ij}d_j &= -t_{ii}^2d_i + ct_{ii}d_i^{-2} - c(q-q^{-1})d_i^{-2},
\tag{b}\\
\sum_{j\ne i,k}v_{ijk} &= -t_{ii} - t_{kk} + (q-q^{-1}),
\quad\text{for $i\ne k$,}\tag{c}\\
\sum_{j\ne i,k}v_{ijk}d_j &= -t_{ii}d_i - t_{kk}d_k + cd_i^{-1}d_k^{-1},
\quad \text{for $i\ne k$}.\tag{d}
\end{align}
\end{prop}

\begin{proof}
Equations (a) and (b) are obtained by comparing the $(i,i)$ entries on
each side of the matrix equations $T^2 = (q-q^{-1})T + q$ and
$TDT = c D^{-1}T^{-1}D^{-1}$. Equations (c) and (d) are obtained by
comparing the $(i,k)$ entries.
\end{proof}

The equations
\begin{equation}\label{eqn:transf}
v_{ijk} = \frac{v_{1ij}v_{1jk}v_{1ki} }{u_{ik} },
\quad
v_{1ji}= \frac{u_{ij}}{v_{1ij}},
\quad
t_{ij} = \frac{t_{1j}v_{1ij}}{t_{1i}},
\end{equation}
imply 
that all of the values in Proposition~\ref{prop:linear}
are determined once we know $u_{ij}$ and $v_{1ij}$ for
$i < j$ and $t_{1i}$ for $1<i\le r$. 

In view of the relations (\ref{eqn:relations}) we may apply
Proposition~\ref{prop:linear} to the matrices $T^k_\lambda$ and
$D^k_\lambda$.  Equations (a) and (b) of Proposition~\ref{prop:linear}
give 
\begin{align*}
&\text{
$4$ equations in the single variable $u_{12}$ when $\dim(T^k_\lambda)=2$,} \\
&\text{
$6$ equations in the $3$ variables $u_{ij}$ when $\dim(T^k_\lambda)=3$,} \\
&\text{
$8$ equations in the $6$ variables $u_{ij}$ when $\dim(T^k_\lambda)=4$,} \\
&\text{
$12$ equations in the $15$ variables $u_{ij}$ when $\dim(T^k_\lambda)=6$.}
\end{align*}
These equations are sufficient to determine the products $u_{ij}=t_{ij}t_{ji}$
for all $T^k_\lambda$ except $T^{25}_{(21)}$ where we have
$\dim(T^{25}_{(21)})=6$. 

If $\dim(T_\lambda^k)=2,3$, or $4$ we use the linear equations in (a) and
(b) of Proposition~\ref{prop:linear} to solve for the $u_{ij}$.
Then we use the equations in (\ref{eqn:transf})
to write the equations in (c) and (d) of Proposition~\ref{prop:linear}
in terms of the variables $t_{1i}$, and $v_{1ij}$, $i<j$.
After doing this we are able to use the subset of the
resulting equations which 
are linear in the $v_{1ij}$ to uniquely determine the values of
the $v_{1ij}$, $i<j$.
This determines the $T^k_\lambda$ up to the choice of the $t_{1i}$.
Finally, the equations resulting from 
the condition $T_3T_4T_3=T_4T_3T_4$ 
force certain relations between the
$T^k_\lambda$ for fixed $k$ and different $\lambda$. 
For each fixed $k$ we picked out a few nice equations resulting from 
this condition to determine the $T^k_\lambda$ completely for all $\lambda$.  
This completely determined the representations $\varphi^k$ for all
$k$ except $k=25$.

The case of $\varphi^{25}$ is slightly more complex.  We used the
same methods as above to determine the matrices $T^{25}_\lambda$
in terms of the variables $t_{1i}$ for each $\lambda$ except
$T^{25}_{(21)}$.  
In the case of $T^{25}_{(21)}$ we have $\dim(T^{25}_{(21)})=6$ and the  
system of 12 equations obtained from (a) and (b)
of Proposition~\ref{prop:linear} 
is a rank 11 system in the 15 unknowns $u_{ij}$.
These linear equations 
can be used to write 11 of the $u_{ij}$ variables in terms of the other 4.  
Next we chose the nicest equations resulting from
(c) and (d) of Proposition~\ref{prop:linear} {\it and} the condition
$T_4T_3T_4=T_3T_4T_3$ and used 
Maple~\cite{char-etal:1991} 
to solve these equations.
These equations are quite nontrivial and we found that we needed to choose
these equations carefully in order to stay within the bounds of
the capability of Maple.
In this way we determined the matrices $T^{25}_\lambda$, for all $\lambda$,
and thus determined $\varphi^{25}$ completely.

%\vfill\pagebreak[4]
\section{Branching rules}\label{secn:branch}

\medskip
\subsection{The branching rules from $HA_2$ to $HA_1$}

\[
\BeginTable
\BeginFormat
| ml | mc | mc | .
\_
|\text{$\varphi^\lambda$} | \text{dim} | \text{Restriction to $HA_1$} |\\+22
\_
|\varphi^{(3)}   |1|\varphi^{(2)}|\\+20
|\varphi^{(2,1)} |2|\varphi^{(2)} \oplus \varphi^{(1^2)}|\\
|\varphi^{(1^3)} |1|\varphi^{(1^2)}|\\+02
\_
\EndTable
\]

\medskip
\subsection{The branching rules from $HB_3$ to $HA_2$.}

\label{tbl:1}
\[
\BeginTable
\BeginFormat
| ml | mc | mc | .
\_
|\text{$\varphi^\mu$} | \text{dim} | \text{Restriction to $HA_2$} |\\+22
\_
|\varphi^{(3),\emptyset}|1|\varphi^{(3)}|\\+20
|\varphi^{(1^3),\emptyset}|1|\varphi^{(1^3)}|\\
|\varphi^{\emptyset,(3)}|1|\varphi^{(3)}|\\
|\varphi^{\emptyset,(1^3)}|1|\varphi^{(1^3)}|\\+02
\_
|\varphi^{(21),\emptyset}|2|\varphi^{(21)}|\\+20
|\varphi^{\emptyset,(21)}|2|\varphi^{(21)}|\\+02
\_
|\varphi^{(2),(1)}|3|\varphi^{(3)}\oplus\varphi^{(21)}|\\+20
|\varphi^{(1^2),(1)}|3|\varphi^{(21)}\oplus\varphi^{(1^3)}|\\
|\varphi^{(1),(2)}|3|\varphi^{(3)}\oplus\varphi^{(21)}|\\
|\varphi^{(1),(1^2)}|3|\varphi^{(21)}\oplus\varphi^{(1^3)}|\\+02
\_
\EndTable
\]

%\vfill\pagebreak[2]

\subsection{The branching rules from $HF_4$ to $HB_3$.}

The 
\vadjust{\nobreak}
bands in this table separate
the orbits of the group of field automorphisms
$\(\alpha_p,\alpha_q\)$, see \ref{secn:f4rep}.
\nobreak
\label{tbl:2}
\[
\BeginTable
\OpenUp{2}{2}
\BeginFormat
|4 ml |4 mc | mc |4 .
\_4
|\text{$\varphi^k$} |\text{dim}| \text{Restriction to $HB_3$} |\\+22
\_
|\varphi^{1}  |1|\varphi^{(3),\emptyset}  |\\
|\varphi^{2}  |1|\varphi^{(1^3),\emptyset}|\\
|\varphi^{3}  |1|\varphi^{\emptyset,(3)}  |\\
|\varphi^{4}  |1|\varphi^{\emptyset,(1^3)}|\\+02
\_
|\varphi^{5}  |2|\varphi^{(21),\emptyset}|\\
|\varphi^{6}  |2|\varphi^{\emptyset,(21)}|\\+02
\_
|\varphi^{7}  |2|\varphi^{(3),\emptyset}\oplus\varphi^{\emptyset,(3)}|\\
|\varphi^{8}  |2|\varphi^{(1^3),\emptyset}\oplus\varphi^{\emptyset,(1^3)}|\\+02
\_
|\varphi^{9}  |4|\varphi^{(21),\emptyset}\oplus\varphi^{\emptyset,(21)}|\\+02
\_
|\varphi^{10} |9|\varphi^{(3),\emptyset} \oplus \varphi^{(21),\emptyset} \oplus
   \varphi^{(2),(1)} \oplus \varphi^{(1),(2)}|\\
|\varphi^{11} |9|\varphi^{(21),\emptyset} \oplus \varphi^{(1^3),\emptyset} \oplus
   \varphi^{(1^2),(1)} \oplus \varphi^{(1),(1^2)}|\\
|\varphi^{12} |9|\varphi^{(2),(1)} \oplus \varphi^{(1),(2)} \oplus
   \varphi^{\emptyset,(3)} \oplus \varphi^{\emptyset,(21)}|\\
|\varphi^{13} |9|\varphi^{(1^2),(1)} \oplus \varphi^{(1),(1^2)} \oplus
   \varphi^{\emptyset,(21)} \oplus \varphi^{\emptyset,(1^3)}|\\+02
\_
|\varphi^{14} |6|\varphi^{(1^2),(1)} \oplus \varphi^{(1),(2)}|\\
|\varphi^{15} |6|\varphi^{(2),(1)} \oplus \varphi^{(1),(1^2)}|\\+02
\_
|\varphi^{16}|12|\varphi^{(2),(1)} \oplus \varphi^{(1^2),(1)} \oplus
                    \varphi^{(1),(2)} \oplus \varphi^{(1),(1^2)}|\\+02
\_
|\varphi^{17} |4|\varphi^{(3),\emptyset}\oplus\varphi^{(2),(1)}|\\
|\varphi^{18} |4|\varphi^{(1^3),\emptyset}\oplus\varphi^{(1^2),(1)}|\\
|\varphi^{19} |4|\varphi^{(1),(2)}\oplus\varphi^{\emptyset,(3)}|\\
|\varphi^{20} |4|\varphi^{(1),(1^2)}\oplus\varphi^{\emptyset,(1^3)}|\\+02
\_
|\varphi^{21} |8|
   \varphi^{(21),\emptyset}\oplus\varphi^{(2),(1)}\oplus\varphi^{(1^2),(1)}|\\
|\varphi^{22} |8|
   \varphi^{(1),(2)}\oplus\varphi^{(1),(1^2)}\oplus\varphi^{\emptyset,(21)}|\\+02
\_
|\varphi^{23} |8|\varphi^{(3),\emptyset}\oplus\varphi^{(2),(1)}\oplus
                 \varphi^{(1),(2)}\oplus\varphi^{\emptyset,(3)}|\\
|\varphi^{24} |8|\varphi^{(1^3),\emptyset}\oplus\varphi^{(1^2),(1)}\oplus
                 \varphi^{(1),(1^2)}\oplus\varphi^{\emptyset,(1^3)}|\\+02
\_
|\varphi^{25}|16|
   \varphi^{(21),\emptyset}\oplus\varphi^{(2),(1)}\oplus\varphi^{(1^2),(1)}\oplus
   \varphi^{(1),(2)}\oplus\varphi^{(1),(1^2)}\oplus\varphi^{\emptyset,(21)}|\\+02
\_4
\EndTable
\]
%\vfill\pagebreak[2]

\section{Seminormal Representations for $HA_1$, $HA_2$,
$HB_3$, and $HF_4$}

\subsection{The Iwahori-Hecke algebra $HA_1$}
\label{secn:a1rep}
%\null\hfil

The irreducible representations $\varphi^\lambda$ of $HA_1$ are
indexed by the partitions $\lambda$ of $2$ and we have
\[
\varphi^{(2)}(T_1) = (p)\quad\text{and}\quad
\varphi^{(1^2)}(T_1) = (-p^{-1}).
%\varphi^\lambda(T_{w_{0,1}}) = p^{\CT(\lambda)}\Id.
\]

\subsection{The Iwahori-Hecke algebra $HA_2$}
\label{secn:a2rep}
%\null\hfil

The irreducible representations $\varphi^\lambda$ of $HA_2$ are
indexed by the partitions of $3$ and can be given explicitly by
\begin{align*}
\varphi^{(3)}(T_1) &= (p),&
\varphi^{(3)}(T_2) &= (p), \\
\varphi^{(21)}(T_1) &= \diag(p,-p^{-1}), &
\varphi^{(21)}(T_2) &= M_2(p,\alpha), \\
\varphi^{(1^3)}(T_1) &= (-p^{-1}),&
\varphi^{(1^3)}(T_2) &= (-p^{-1}), 
\end{align*}
where
\[
M_2(p,\alpha) = -\frac{1}{[2]_p}\begin{pmatrix}
p^{-2}  &\alpha([2]_p-1)\\[4pt]
\frac{1}{\alpha}([2]_p+1)  &-p^2
\end{pmatrix}
\]
and $[2]_p = p+p^{-1}$.  The variable $\alpha$ is a free parameter,
see Lemma~\ref{prop:param}.

\subsection{The Iwahori-Hecke algebra $HB_3$}
\label{secn:b3rep}

The irreducible representations $\varphi^\mu=\varphi^{\alpha,\beta}$ of
$HB_3(p^2,q^2)$ are indexed by pairs of partitions $\mu=(\alpha,\beta)$ such
that $|\alpha|+|\beta|=3$. 
Let $\diag(A,B,\ldots,C)$ denote the block
diagonal matrix with the matrices $A, B,\ldots, C$ in order on the diagonal.
Then, using the notation
\[
[2]_x = x+x^{-1},
\qquad
[3]_x = x^2+1+x^{-2},
\qquad\text{and}\qquad
[0]_x=x-x^{-1},
\]
irreducible seminormal representations of the Iwahori-Hecke algebra
$HB_3$ can be given explicitly as follows:
\begin{align*}
\varphi^{(3),\emptyset}(T_1) &= (p),&
\varphi^{(3),\emptyset}(T_2) &= (p),&
\varphi^{(3),\emptyset}(T_3) &= (q), \\
\varphi^{(1^3),\emptyset}(T_1) &= (-p^{-1}),&
\varphi^{(1^3),\emptyset}(T_2) &= (-p^{-1}),&
\varphi^{(1^3),\emptyset}(T_3) &= (q), \\
\varphi^{\emptyset, (3)}(T_1) &= (p),&
\varphi^{\emptyset, (3)}(T_2) &= (p),&
\varphi^{\emptyset, (3)}(T_3) &= (-q^{-1}), \\
\varphi^{\emptyset, (1^3)}(T_1) &= (-p^{-1}),&
\varphi^{\emptyset, (1^3)}(T_2) &= (-p^{-1}),&
\varphi^{\emptyset, (1^3)}(T_3) &= (-q^{-1}), \\
\varphi^{(21),\emptyset}(T_1) &= \diag(p,-p^{-1}),&
\varphi^{(21),\emptyset}(T_2) &=  M_2(p,1),&
\varphi^{(21),\emptyset}(T_3) &=  \diag(q,q), \\
\varphi^{\emptyset, (21)}(T_1) &= \diag(p,-p^{-1}),&
\varphi^{\emptyset, (21)}(T_2) &=  M_2(p,1),&
\varphi^{\emptyset, (21)}(T_3) &=  \diag(-q^{-1},-q^{-1}).
\end{align*}

\begin{align*}
\varphi^{(2),(1)}(T_1) &= \diag(p,p,-p^{-1}),&
\varphi^{(2),(1)}(T_2) &= \diag(p,M_2(p,1)), \\
\varphi^{(1^2),(1)}(T_1) &= \diag(p,-p^{-1},-p^{-1}), &
\varphi^{(1^2),(1)}(T_2) &= \diag(M_2(p,1),-p^{-1}), \\
\varphi^{(1),(2)}(T_1) &= \diag(p,p,-p^{-1}), &
\varphi^{(1),(2)}(T_2) &= \diag(p,M_2(p,1)), \\
\varphi^{(1),(1^2)}(T_1) &= \diag(p,-p^{-1},-p^{-1}), &
\varphi^{(1),(1^2)}(T_2) &= \diag(M_2(p,1),-p^{-1}). 
\end{align*}

\begin{align*}
\varphi^{(2),(1)}(T_3) &= \diag(M_{(2),(1)},q), \\
\varphi^{(1^2),(1)}(T_3) &= \diag(q,M_{(1^2),(1)}), \\
\varphi^{(1),(2)}(T_3) &= \diag(M_{(1),(2)},-q^{-1}), \\
\varphi^{(1),(1^2)}(T_3) &= \diag(-q^{-1},M_{(1),(1^2)}), 
\end{align*}
where 
\begin{align*}
M_{(2),(1)} &= \frac{1}{[3]_p}
\begin{pmatrix}
q+p^{-2}[0]_q& -[2]_p[2]_{p/q}\\[6pt]
-[2]_{p^2q}&-q^{-1}+p^2[0]_q
\end{pmatrix}, \\
\\
M_{(1^2),(1)} &= \frac{1}{[3]_p}
\begin{pmatrix}
-q^{-1}+p^{-2}[0]_q& -[2]_{p^2/q}\\[6pt]
-[2]_p[2]_{p^2q}&q+p^2[0]_q
\end{pmatrix}, \\
\\
M_{(1),(2)} &= \frac{1}{[3]_p}
\begin{pmatrix}
-q^{-1}+p^{-2}[0]_q& [2]_p[2]_{pq}\\[6pt]
% Note that if the next character were [ it would be interpreted by
% LaTeX as the beginning of a command to increase the line spacing,
% hence the braces around [2]_ ..  But now I really do have [4pt],
% so this precaution is actually not necessary!
{[2]_{p^2/q}}&q+p^2[0]_q
\end{pmatrix}, \\
\\
M_{(1),(1^2)} &= \frac{1}{[3]_p}
\begin{pmatrix}
q+p^{-2}[0]_q& [2]_{p^2q}\\[6pt]
{[2]_p[2]_{p/q}}&-q^{-1}+p^2[0]_q
\end{pmatrix}.
\end{align*}

\subsection{The Iwahori-Hecke algebra $HF_4$}
\label{secn:f4rep}

Let $\alpha_p$ be the automorphism of $\Q(p,q)$ which fixes $q$ and
sends $p$ to $-p^{-1}$. Similarly, let $\alpha_q$ be the automorphism
which fixes $p$ and sends $q$ to $-q^{-1}$. These field automorphisms
act on the entries of the matrices $\varphi^\lambda(T_i)$ and thereby
permute the representations $\varphi^\lambda$. 
The representation resulting from the application of a field
automorphism to a representation in seminormal form
may no longer be seminormal. In order to bring the
representation back to seminormal form it may be necessary to conjugate
by a permutation matrix $\pi$. 
The orbits of the irreducible representations under the action
of $\alpha_p$ and $\alpha_q$ and the permutations $\pi$ for
conjugating back to seminormal form are given in the 
following table.  If $\varphi$ is a representation of $HF_4$
and $\pi$ is a permutation then we shall let $\pi\circ \varphi$ denote the
representation determined
by $(\pi\circ \varphi)(h)=\pi\varphi(h)\pi^{-1}$,
for all $h\in HF_4$.

\label{tbl:3}
\[
\BeginTable
\OpenUp{2}{2}
\BeginFormat
| ml | mc |.
\_
| \varphi^k |\text{Orbit of $\(\alpha_p,\alpha_q\)$}|\\
\_
|\varphi^{1}  |
\varphi^2 = \alpha_p\varphi^1,\quad
\varphi^3 = \alpha_q\varphi^1,\quad
\text{and}\quad
\varphi^4 = \alpha_p\alpha_q\varphi^1|\\
\_
|\varphi^{5}  |\varphi^6 = \alpha_q\varphi^5|\\
\_
|\varphi^{7}  |\varphi^8 = \alpha_p\varphi^7|\\
\_
|\varphi^{10} |\varphi^{11}=\pi_{11}\circ(\alpha_p\varphi^{10}),
\quad \text{where $\pi_{11} = (1,3)(4,6)(7,9)$}|\\
| |\varphi^{12} = \pi_{12}\circ(\alpha_q\varphi^{10}),
\quad \text{where $\pi_{12} = (1,7)(2,8)(3,9)$}|\\
| |\varphi^{13} = \pi_{13}\circ(\alpha_p\alpha_q\varphi^{10}),
\quad \text{where $\pi_{13} = (1,9)(2,8)(3,7)(4,6)$}|\\
\_
|\varphi^{14} |\varphi^{15} = \pi_{15}\circ(\alpha_p\varphi^{14}),
\quad \text{where $\pi_{15} = (1,3)(4,6)$}|\\
\_
|\varphi^{17} |\varphi^{18} = \pi_{18}\circ(\alpha_p\varphi^{17}),
\quad \text{where $\pi_{18} = (2,4)$}|\\
| |\varphi^{19} = \pi_{19}\circ(\alpha_q\varphi^{17}),
\quad \text{where $\pi_{19} = (1,4,3,2)$}|\\
| |\varphi^{20} = \pi_{20}\circ(\alpha_p\alpha_q\varphi^{17}),
\quad \text{where $\pi_{20} = (1,4)(2,3)$}|\\
\_
|\varphi^{21} |\varphi^{22} = \pi_{22}\circ(\alpha_q\varphi^{21}),
\quad \text{where $\pi_{22} = (1,7,5,3)(2,8,6,4)$}|\\
\_
|\varphi^{23} |\varphi^{24} = \pi_{24}\circ(\alpha_p\varphi^{23}),
\quad \text{where $\pi_{24} = (2,4)(5,7)$}|\\
\_
\EndTable
\]

Let $w_{0,1}$, $w_{0,2}$, $w_{0,3}$ and $w_{0,4}$
be the longest elements in the Weyl groups $WA_1$, $WA_2$,
$WB_3$ and $WF_4$, respectively. Let
\begin{align*}
D_1 &= T_{w_{0,1}} = T_1, \\
D_2 &= T_{w_{0,2}}^2 = (T_1T_2T_1)^2, \\
D_3 &= T_{w_{0,3}} = (T_3T_2T_1)^3, \\
D_4 &= T_{w_{0,4}} = (T_4T_{w_{0,3}})^3 T_{w_{0,2}}^{-2} 
\end{align*}
in $HF_4$.
The following tables give the values of $\varphi^k(D_j)$,
for one representative from each equivalence class of representations.
The rest of the matrices $\varphi^k(D_j)$ are easily obtained
by applying the automorphisms $\alpha_p$ and $\alpha_q$ and
conjugating by a permutation $\pi$ as indicated in \ref{tbl:3}
above.

\[
\begin{aligned}
  \varphi^{1}(T_1) &= (p),\\
  \varphi^{1}(T_{w_{0,2}}^2) &= (p^6),\\
  \varphi^{1}(T_{w_{0,3}}) &= (p^6q^3),\\
  \varphi^{1}(T_{w_{0,4}}) &= (p^{12}q^{12} ),
\end{aligned}
\qquad\qquad\qquad\qquad
\begin{aligned}
  \varphi^{5}(T_1) &= \diag(p,-p^{-1}),\\
  \varphi^{5}(T_{w_{0,2}}^2) &= \Id,\\
  \varphi^{5}(T_{w_{0,3}}) &= q^3\,\Id,\\
  \varphi^{5}(T_{w_{0,4}}) &= q^{12}\,\Id,
\end{aligned}\qquad
\]

\[
\begin{aligned}
  \varphi^{7}(T_1) &= p\,\Id,\\
  \varphi^{7}(T_{w_{0,2}}^2) &= p^6\,\Id,\\
  \varphi^{7}(T_{w_{0,3}}) &= \diag(p^6q^3,-p^6q^{-3}),\\
  \varphi^{7}(T_{w_{0,4}}) &= p^{12}\,\Id,
\end{aligned}
\qquad
\begin{aligned}
  \varphi^{9}(T_1) &= \diag(p,-p^{-1},p,-p^{-1}),\\
  \varphi^{9}(T_{w_{0,2}}^2) &= \Id,\\
  \varphi^{9}(T_{w_{0,3}}) &= \diag(q^3,q^3,-q^{-3},-q^{-3}),\\
  \varphi^{9}(T_{w_{0,4}}) &= \Id,
\end{aligned}
\]

\begin{align*}
  \varphi^{10}(T_1) &= \diag(p,p,-p^{-1},p,p,-p^{-1},p,p,-p^{-1},),
  \hfil \\
  \varphi^{10}(T_{w_{0,2}}^2) &= \diag(p^6,1,1,p^6,1,1,p^6,1,1),\\
  \varphi^{10}(T_{w_{0,3}}) &= \diag(p^6q^3,q^3,q^3,-p^2q,-p^2q,-p^2q,
          p^2q^{-1},p^2q^{-1},p^2q^{-1}),\\
  \varphi^{10}(T_{w_{0,4}}) &= p^4q^4\,\Id,
\end{align*}

\begin{align*}
  \varphi^{14}(T_1) &= \diag(p,-p^{-1},-p^{-1},p,p,-p^{-1}),
  \hfil \\
  \varphi^{14}(T_{w_{0,2}}^2) &= \diag(1,1,p^{-6},p^6,1,1),\\
  \varphi^{14}(T_{w_{0,3}}) &= \diag(-p^{-2}q,-p^{-2}q,-p^{-2}q,
        p^2q^{-1},p^2q^{-1},p^2q^{-1}),\\
  \varphi^{14}(T_{w_{0,4}}) &= \Id,
\end{align*}

\begin{align*}
  \varphi^{16}(T_1) &= \diag(p,p,-p^{-1},p,-p^{-1},-p^{-1},p,p,-p^{-1},
         p,-p^{-1},-p^{-1}),
	 \hfil \\
  \varphi^{16}(T_{w_{0,2}}^2) &= \diag(p^6,1,1,1,1,p^{-6},p^6,
         1,1,1,1,p^{-6}),\\
  \varphi^{16}(T_{w_{0,3}}) &= \diag(-p^2q,-p^2q,-p^2q,
          -p^{-2}q,-p^{-2}q,-p^{-2}q,p^2q^{-1},p^2q^{-1},p^2q^{-1},\\
          &\qquad p^{-2}q^{-1},p^{-2}q^{-1},p^{-2}q^{-1}),\\
  \varphi^{16}(T_{w_{0,4}}) &= \Id,
\end{align*}

\begin{align*}
  \varphi^{17}(T_1) &= \diag(p,p,p,-p^{-1}),
  \hfil \\
  \varphi^{17}(T_{w_{0,2}}^2) &= \diag(p^6,p^6,1,1),\\
  \varphi^{17}(T_{w_{0,3}}) &= \diag(p^6q^3,-p^2q,-p^2q,-p^2q),\\
  \varphi^{17}(T_{w_{0,4}}) &= -p^6q^6\,\Id,
\end{align*}

\begin{align*}
  \varphi^{21}(T_1) &= \diag(p,-p^{-1},p,p,-p^{-1},p,-p^{-1},-p^{-1}),
  \hfil \\
  \varphi^{21}(T_{w_{0,2}}^2) &= \diag(1,1,p^6,1,1,1,1,p^{-6}),\\
  \varphi^{21}(T_{w_{0,3}}) &= \diag(q^3,q^3,-p^2q,-p^2q,-p^2q,-p^{-2}q,
          -p^{-2}q,-p^{-2}q),\\
  \varphi^{21}(T_{w_{0,4}}) &= -q^6\,\Id,
\end{align*}

\begin{align*}
  \varphi^{23}(T_1) &= \diag(p,p,p,-p^{-1},p,p,-p^{-1},p),
  \hfil \\
  \varphi^{23}(T_{w_{0,2}}^2) &= \diag(p^6,p^6,1,1,p^6,1,1,p^6),\\
  \varphi^{23}(T_{w_{0,3}}) &= \diag(p^6q^3,-p^2q,-p^2q,-p^2q,
             p^2q^{-1},p^2q^{-1},p^2q^{-1},-p^6q^{-3}),\\
  \varphi^{23}(T_{w_{0,4}}) &= -p^6\,\Id,
\end{align*}

\begin{align*}
\varphi^{25}(T_1) &= \diag(p,-p^{-1},p,p,-p^{-1},p,-p^{-1},-p^{-1},
        p,p,-p^{-1},p,-p^{-1},-p^{-1},p,-p^{-1}), \hfil \\
\varphi^{25}(T_{w_{0,2}}^2) &= \diag(1,1,p^6,1,1,1,1,p^{-6},p^6,
        1,1,1,1,p^{-6},1,1),\\
\varphi^{25}(T_{w_{0,3}}) &= \diag(q^3,q^3,-p^2q,-p^2q,-p^2q,
        -p^{-2}q,-p^{-2}q,-p^{-2}q,p^2q^{-1},p^2q^{-1},p^2q^{-1},\\
        &\qquad p^{-2}q^{-1},p^{-2}q^{-1},p^{-2}q^{-1},-q^{-3},q^{-3}),\\
\varphi^{25}(T_{w_{0,4}}) &= -\Id. 
\end{align*}

Using the methods described in the previous sections we
have produced matrices $\varphi^k(T_i)$ giving 
the 25 irreducible representations $\varphi^k$ of $HF_4$.
The following tables give the values of $\varphi^k(T_i)$,
for one representative from each equivalence class of representations.
The rest of the matrices $\varphi^k(T_i)$ are obtained
by applying the automorphisms $\alpha_p$ and $\alpha_q$ and
conjugating by a permutation $\pi$ as indicated in \ref{tbl:3}
above.

We shall use the notations
\[
[2]_x = x+x^{-1},
\qquad
[3]_x = x^2+1+x^{-2},
\qquad\text{and}\qquad
[0]_x=x-x^{-1},
\]
and the notation
\[
\varphi^k(T_i)^{[a_1,a_2,\ldots,a_r]}
\]
will denote the $r\times r$ submatrix of $\varphi^k(T_i)$ which
is formed by the intersection of the $a_1,\ldots, a_r$th
rows and columns.  
The notation $\diag(A,B,\ldots,C)$ will denote the block
diagonal matrix with the matrices $A, B,\ldots, C$ in order 
along the diagonal.
The matrix $M_2(x,y)$ will be as given in \ref{secn:a2rep},
the matrices $M_{\alpha,\beta}$ are as
given in \ref{secn:b3rep} and the variables
$\alpha$, $\beta$, $\xi$, $\theta$, and $\eta$ 
are free parameters, see Lemma~\ref{prop:param}.
{\it Any entries of the matrices $\varphi^k(T_i)$
which are not given explicitly below are taken to be $0$.}

\filbreak
\subsubsection*{The representations $\varphi^1$ and $\varphi^5$}
\[
\begin{aligned}
\varphi^1(T_1) &= (p),\\
\varphi^1(T_2) &= (p),\\
\varphi^1(T_3) &= (q),\\
\varphi^1(T_4) &= (q),
\end{aligned}
\qquad\qquad\qquad
\begin{aligned}
\varphi^5(T_1) &= \diag(p,-p^{-1}),\\
\varphi^5(T_2) &= M_2(p,1),\\
\varphi^5(T_3) &= \diag(q,q),\\
\varphi^5(T_4) &= \diag(q,q).
\end{aligned}\qquad
\]

\filbreak
\subsubsection*{The representations $\varphi^7$ and $\varphi^{9}$}
\[
\begin{aligned}
\varphi^7(T_1) &= \diag(p,p),\\
\varphi^7(T_2) &= \diag(p,p),\\
\varphi^7(T_3) &= \diag(q,-q^{-1}),\\
\varphi^7(T_4) &= M_2(q,\alpha),
\end{aligned}
\qquad
\begin{aligned}
\varphi^9(T_1) &= \diag(p,-p^{-1},p,-p^{-1}),\\
\varphi^9(T_2) &= \diag(M_2(p,1),M_2(p,1)),\\
\varphi^9(T_3) &= \diag(q,q,-q^{-1},q^{-1}),\\
\varphi^9(T_4)^{[1,3]} &= \varphi^9(T_4)^{[2,4]} = M_2(q,\alpha).
\end{aligned}
\]

\filbreak
\subsubsection*{The representation $\varphi^{10}$}
\begin{align*}
\varphi^{10}(T_1) &= \diag(p,p,-p^{-1},p,p,-p^{-1},p,p,-p^{-1}),\\
\varphi^{10}(T_2) &= \diag(p,M_2(p,1),p,M_2(p,1),p,M_2(p,1)),\\
\varphi^{10}(T_3) &= \diag(q,q,q,M_{(2,1)},q,M_{(1,2)},q),\\
\varphi^{10}(T_4)^{[1,4,7]} &= M_{10},\\
\varphi^{10}(T_4)^{[2,5,8]} &= \varphi^{10}(T_4)^{[3,6,9]} = N_{10},
\end{align*}

where
\[
M_{10} = \dfrac{1}{[2]_q[2]_{p^2q}}\begin{pmatrix}
    p^{-2}q^{-1}[2]_q[0]_q&
   -[2]_q[2]_{p^2q^2}\xi\eta^{-1}&-[2]_q[2]_{p^2q^2}\xi\\[8pt]
%\noalign{\vskip\alignskip}
   -[2]_{p^2q^{-1}}\eta\xi^{-1}&
    [2]_{p^2q}+p^2q[2]_{q}[0]_q&-[2]_{p^2q^{-1}}\eta\\[8pt]
%\noalign{\vskip\alignskip}
   -[2]_{p^2q}\xi^{-1}& -[2]_{p^2q}\eta^{-1}&q^2[2]_{p^2q}
\end{pmatrix}
\]

and
\[
N_{10} = \dfrac{1}{[2]_q[2]_{pq^{-1}}}\begin{pmatrix}
    pq^{-1}[2]_q[0]_q&
   -[2]_q[2]_{pq^{-2}}\theta\eta^{-1}&-[2]_q[2]_{pq^{-2}}\theta\\[8pt]
%\noalign{\vskip\alignskip}
   -[2]_{pq}\eta\theta^{-1}&
    [2]_{pq^{-1}}+p^{-1}q[2]_{q}[0]_q&-[2]_{pq}\eta\\[8pt]
%\noalign{\vskip\alignskip}
   -[2]_{pq^{-1}}\theta^{-1}&-[2]_{pq^{-1}}\theta^{-1}&q^2[2]_{pq^{-1}}
\end{pmatrix}.
\]

\filbreak
\subsubsection*{The representation $\varphi^{14}$}
\begin{align*}
\varphi^{14}(T_1) &= \diag(p,-p^{-1},-p^{-1},p,p,-p^{-1}),\\
\varphi^{14}(T_2) &= \diag(M_2(p,1),-p^{-1},p,M_2(p,1)),\\
\varphi^{14}(T_3) &= \diag(q,M_{(1^2,1)},M_{(1,2)},-q^{-1}),\\
\varphi^{14}(T_4)^{[3]} &= -q^{-1},\\
\varphi^{14}(T_4)^{[4]} &= q,\\
\varphi^{14}(T_4)^{[1,5]} &= M_{14},\\
\varphi^{14}(T_4)^{[2,6]} &= M_{14},
\end{align*}

where
\[
M_{14} = \dfrac{1}{[2]_{p^2q^{-1}}}
\begin{pmatrix}
1+p^2q^{-1}[0]_q&-[3]_p\alpha\\[6pt]
%\noalign{\vskip\alignskip}
(1-[2]_{p^2q^{-2}})\alpha^{-1}&-1+p^{-2}q[0]_q
\end{pmatrix}.
\]

\filbreak
\subsubsection*{The representation $\varphi^{16}$}
\begin{align*}
\varphi^{16}(T_1) &=
\diag(p,p,-p^{-1},p,-p^{-1},-p^{-1},p,p,-p^{-1},p,-p^{-1},-p^{-1}),\\
\varphi^{16}(T_2) &= \diag(p,M_2(p,1),M_2(p,1),-p^{-1},p,M_2(p,1),
           M_2(p,1)-p^{-1}),\\
\varphi^{16}(T_3) &=
\diag(M_{(2,1)},q,q,M_{(1^2,1)},M_{(1,2)},-q^{-1},-q^{-1},M_{(1,1^2)}),\\
\varphi^{16}(T_4)^{[1,7]} &= M_{16}(\xi),\\
\varphi^{16}(T_4)^{[6,12]} &= M_{16}(\eta),\\
\varphi^{16}(T_4)^{[2,4,8,10]} &= \varphi^{16}(T_4)^{[3,5,9,12]} = N_{16},
\end{align*}

where
\[
M_{16}(\alpha) = \dfrac{1}{[2]_q}
\begin{pmatrix}
1+q^{-1}[0]_q& -3\alpha\\[6pt]
%\noalign{\vskip\alignskip}
-[3]_{q^2}/\alpha [3]_q&-1+q[0]_q
\end{pmatrix},
\]

\[
N_{16} = \dfrac{1}{[2]_p[2]_q}\begin{pmatrix}
 f_{16}(p,q)&\dfrac{3[2]_{pq}\xi\theta}{[2]_{p^2/q}\eta}&
 \dfrac{3[2]_{pq}\xi}{[2]_{p^2/q}}&\dfrac{3[2]_{pq}\xi\theta}{[2]_{p^2q}}\\[14pt]
%\noalign{\vskip\alignskip}
 \dfrac{[3]_{p^2}[2]_{p/q}\eta}{[2]_{p^2q}\xi\theta}&-f_{16}(-p^{-1},q)&
 \dfrac{[3]_{p^2}[2]_{p/q}\eta}{[2]_{p^2/q}\theta}&
 \dfrac{3[2]_{p/q}\eta}{[2]_{p^2q}}\\[14pt]
%\noalign{\vskip\alignskip}
 \dfrac{[3]_{q^2}[2]_{p/q}}{[3]_q[2]_{p^2q}\xi}&
 \dfrac{[3]_{q^2}[2]_{p/q}\theta}{[3]_q[2]_{p^2/q}\eta}&
-f_{16}(p,-q^{-1})&-\dfrac{3[2]_{p/q}\theta}{[2]_{p^2q}}\\[14pt]
%\noalign{\vskip\alignskip}
 \dfrac{[3]_{p^2}[3]_{q^2}[2]_{pq}}{3[3]_q[2]_{p^2q}\xi\theta}&
 \dfrac{[3]_{q^2}[2]_{pq}}{[3]_q[2]_{p^2/q}\eta}&
-\dfrac{[3]_{p^2}[2]_{pq}}{[2]_{p^2/q}\theta}&f_{16}(-p^{-1},-q^{-1})
\end{pmatrix},
\]

and
\[
f_{16}(x,y) = \dfrac{-2x/y+xy+1/xy-1/xy^3-y/x-1/x^3y^3+y/x^3}{[2]_{x^2y}}.
\]

\filbreak
\subsubsection*{The representation $\varphi^{17}$}
\begin{align*}
\varphi^{17}(T_1) &= \diag(p,p,p,-p^{-1}),\\
\varphi^{17}(T_2) &= \diag(p,p,M_2(p,1)),\\
\varphi^{17}(T_3) &= \diag(q,M_{(2,1)},q),\\
\varphi^{17}(T_4) &= \diag(M_{17},q,q),
\end{align*}

where
\[
M_{17} = \dfrac{1}{[2]_{p^2q}}\begin{pmatrix}
1+p^{-2}q^{-1}[0]_q&-[3]_p\alpha\\
\noalign{\vskip\alignskip}
(1-[2]_{p^2q^2})\alpha^{-1}&-1 + p^2q[0]_q
\end{pmatrix}.
\]

\filbreak
\subsubsection*{The representation $\varphi^{21}$}
\begin{align*}
\varphi^{21}(T_1) &= \diag(p,-p^{-1},p,p,-p^{-1},p,-p^{-1},-p^{-1}),\\
\varphi^{21}(T_2) &= \diag(M_2(p,1),p,M_2(p,1),M_2(p,1),-p^{-1}),\\
\varphi^{21}(T_3) &= \diag(q,q,M_{(2,1)},q,q,M_{(1^2,1)}),\\
\varphi^{21}(T_4)^{[3]} &= \varphi^{21}(T_4)^{[8]} = q,\\
\varphi^{21}(T_4)^{[1,4,6]} &= \varphi^{21}(T_4)^{[2,5,7]} = M_{21},
\end{align*}

where
\[
M_{21} = \dfrac{1}{[2]_p[2]_{pq}[2]_{p/q}}\begin{pmatrix}
(q[2]_{p^2}+q^{-2}[0]_q)[2]_p&-[2]_{q^3}[2]_p\xi\eta^{-1}&-[2]_{q^3}[2]_p\xi\\
[8pt]%\noalign{\vskip\alignskip}
-[3]_p[2]_{pq}\eta\xi^{-1}&(p^{-1}q[2]_p[0]_q+1)[2]_{pq}&-[3]_p[2]_{pq}\eta\\
[8pt]%\noalign{\vskip\alignskip}
-[3]_p[2]_{p/q}\xi^{-1}&-[3]_p[2]_{p/q}\eta^{-1}&(pq[2]_p[0]_q+1)[2]_{p/q}\\
\end{pmatrix}.
\]

\filbreak
\subsubsection*{The representation $\varphi^{23}$}
\begin{align*}
\varphi^{23}(T_1) &= \diag(p,p,p,-p^{-1},p,p,-p^{-1},p),\\
\varphi^{23}(T_2) &= \diag(p,p,M_2(p,1),p,M_2(p,1),p),\\
\varphi^{23}(T_3) &= \diag(q,M_{(2,1)},q,M_{(1,2)},-q^{-1},-q^{-1}),\\
\varphi^{23}(T_4)^{[3,6]} &= \varphi^{23}(T_4)^{[4,7]} =
  M_2(q,\eta/([2]_q-1)\theta),\\
\varphi^{23}(T_4)^{[1,2,5,8]} &= M_{23},
\end{align*}

where
\[
M_{23} = \dfrac{1}{[2]_q}
\begin{pmatrix}
  f_{23}(p,q)&
  \dfrac{[2]_{pq^2}\xi}{[2]_{p^2q}[2]_{pq}\eta}&
  \dfrac{[2]_{p/q}[2]_{pq^2}\xi}{[2]_{p^2q}[2]_{pq}^2\theta}&
  \dfrac{[2]_{pq^2}\xi}{[2]_{p^2q}[2]_{pq}}\\[14pt]
%\noalign{\vskip\alignskip}
  \dfrac{[3]_{p^2}[2]_{p}\eta}{[2]_{p^2q}[2]_{p/q}\xi}&
  g_{23}(p,q)&
  \dfrac{([2]_{q^2}-1)[2]_{p}\eta}{[2]_{pq}[2]_{p^2q}\theta}&
  \dfrac{-[3]_{p^2}[2]_{p}\eta}{[2]_{p/q}[2]_{p^2q}}\\[14pt]
%\noalign{\vskip\alignskip}
  \dfrac{[3]_{q}[3]_{p^2}[2]_{p}\theta}{[2]_{p/q}[2]_{p^2/q}\xi}&
  \dfrac{([2]_{q^2}-1)[3]_{q}[2]_{p}\theta}{[2]_{p/q}[2]_{p^2/q}\eta}&
  -g_{23}(p,-1/q)&
  \dfrac{[3]_{p^2}[2]_{p}\theta}{[2]_{p/q}[2]_{p^2/q}}\\[14pt]
%\noalign{\vskip\alignskip}
  \dfrac{[3]_{q}[2]_{p/q^2}}{[2]_{p/q}[2]_{p^2/q}\xi}&
  \dfrac{-[3]_{q}[2]_{p/q^2}}{[2]_{p/q}[2]_{p^2/q}\eta}&
  \dfrac{[2]_{p/q^2}}{[2]_{pq}[2]_{p^2/q}\theta}&
  -f_{23}(p,-1/q)
\end{pmatrix},
\]

where
\begin{align*}
 f_{23}(x,y) &= \dfrac{y^4 - x^4 y^2 - x^2 y^2 -1}{x^3 y^4 [2]_{xy}[2]_{x^2y}}
\quad\text{and}\\[8pt]
%\noalign{\vskip\alignskip}
 g_{23}(x,y) &= \dfrac{x^4y^6 - x^4y^2 +x^2y^6 - x^2y^4 + x^2y^2 + y^6 + y^2 - 1}{xy^4[2]_{x/y}[2]_{x^2y}}.
\end{align*}

\subsubsection*{The representation $\varphi^{25}$}
\begin{align*}
\varphi^{25}(T_1) &=
\diag(p,-p^{-1},p,p,-p^{-1},p,-p^{-1},-p^{-1},\\
   &\qquad\qquad\qquad\qquad p,p,-p^{-1},p,-p^{-1},-p^{-1},p,-p^{-1}),\\
\varphi^{25}(T_2) &= \diag(M_2(p,1),p,M_2(p,1),M_2(p,1),-p^{-1},\\
   &\qquad\qquad\qquad\qquad p,M_2(p,1),M_2(p,1),-p^{-1},M_2(p,1)),\\
\varphi^{25}(T_3) &= \diag(q,q,M_{(2,1)},q,q,M_{(1^2,1)},
   M_{(1,2)},-q^{-1},-q^{-1},M_{(1,1^2)},-q^{-1},-q^{-1}),\\
\varphi^{25}(T_4)^{[3,9]} &= M_2(q,\alpha/([2]_q-1)\eta),\\
\varphi^{25}(T_4)^{[8,14]} &= M_2(q,\beta/([2]_q-1)\theta),\\
\varphi^{25}(T_4)^{[1,4,6,10,12,15]} &=
\varphi^{25}(T_4)^{[2,5,7,11,13,16]} = M_{25},
\end{align*}

\filbreak
where 
\begin{multline*}
M_{25} = \dfrac{1}{[2]_q}
\left(\begin{matrix}
f_{25}(p,q)&
\dfrac{-[2]_{p^2/q} \xi}{[2]_{p q} [2]_{p/q} \alpha}&
\dfrac{-[2]_{p^2 q} \xi}{[2]_{p q} [2]_{p/q} \beta}\\[14pt] %1
%\noalign{\vskip\alignskip}
\dfrac{-2 [3]_{p} [2]_{p^2/q^2} \alpha}{[2]_{p} [2]_{p/q}[2]_{p^2 q}[2]_{p^2/q}\xi}&
g_{25}(p,q)&
\dfrac{[3]_{p} [2]_{p^2/q^2} \alpha}{[2]_{p} [2]_{p/q} [2]_{p^2/q}\beta}\\[14pt] %2
%\noalign{\vskip\alignskip}
\dfrac{-2 [3]_{p} [2]_{p^2 q^2} \beta}{[2]_{p} [2]_{p q} [2]_{p^2 q} [2]_{p^2/q} \xi}&
\dfrac{[3]_{p} [2]_{p^2 q^2} \beta}{[2]_{p} [2]_{p q} [2]_{p^2 q} \alpha}&
g_{25}(-1/p,q)\\[14pt] %3
%\noalign{\vskip\alignskip}
\dfrac{-2 [3]_{p} [2]_{p^2/q^2} [3]_{q} \eta}{[2]_{p} [2]_{p q} [2]_{p^2 q} [2]_{p^2/q} \xi}&
\dfrac{-([3]_{p}-[3]_{q}+2) [3]_{q} \eta}{[2]_{p} [2]_{p q} [2]_{p^2 q} \alpha}&
\dfrac{[3]_{p} [2]_{p^2/q^2} [3]_{q} \eta}{[2]_{p} [2]_{p q}[2]_{p^2/q} \beta}\\[14pt] %4
%\noalign{\vskip\alignskip}
\dfrac{-2 [3]_{p} [2]_{p^2 q^2} [3]_{q} \theta}{[2]_{p} [2]_{p/q} [2]_{p^2 q} [2]_{p^2/q} \xi}&
\dfrac{[3]_{p} [2]_{p^2 q^2} [3]_{q} \theta}{[2]_{p} [2]_{p/q} [2]_{p^2 q} \alpha}&
\dfrac{-([3]_{p}-[3]_{q}+2) [3]_{q} \theta}{[2]_{p} [2]_{p/q} [2]_{p^2/q} \beta}\\[14pt] %5
%\noalign{\vskip\alignskip}
\dfrac{[2]_{p^2 q^2} [2]_{p^2/q^2}}{[2]_{p q} [2]_{p/q} [2]_{p^2 q} [2]_{p^2/q} \xi}&
\dfrac{[2]_{p^2 q^2}}{[2]_{p q} [2]_{p/q} [2]_{p^2 q} \alpha}&
\dfrac{[2]_{p^2/q^2}}{[2]_{p q} [2]_{p/q} [2]_{p^2/q} \beta}\\ %6
\end{matrix}\right.
\\ \\ \\
\left.\begin{matrix}
\dfrac{-[2]_{p^2 q} \xi}{[2]_{p q} [2]_{p/q} \eta}&
\dfrac{-[2]_{p^2/q} \xi}{[2]_{p q} [2]_{p/q} \theta}&
\dfrac{[3]_{q} [2]_{p^2 q} [2]_{p^2/q} \xi}{[2]_{p q} [2]_{p/q}}\\[14pt]
%\noalign{\vskip\alignskip}
\dfrac{-([3]_{p}-[3]_{q}+2) \alpha}{[2]_{p} [2]_{p/q} [2]_{p^2/q} \eta}&
\dfrac{[3]_{p} [2]_{p^2/q^2} \alpha}{[2]_{p} [2]_{p/q} [2]_{p^2 q} \theta}&
\dfrac{2 [3]_{p} [3]_{q} \alpha}{[2]_{p} [2]_{p/q}}\\[14pt]
%\noalign{\vskip\alignskip}
\dfrac{[3]_{p} [2]_{p^2 q^2} \beta}{[2]_{p} [2]_{p q} [2]_{p^2/q} \eta}&
\dfrac{-([3]_{p}-[3]_{q}+2) \beta}{[2]_{p} [2]_{p q} [2]_{p^2 q} \theta}&
\dfrac{2 [3]_{p} [3]_{q} \beta}{[2]_{p} [2]_{p q}}\\[14pt]
%\noalign{\vskip\alignskip}
-g_{25}(p,-1/q)&
\dfrac{-[3]_{p} [2]_{p^2/q^2} \eta}{[2]_{p} [2]_{p q} [2]_{p^2 q} \theta}&
\dfrac{-2 [3]_{p} [3]_{q} \eta}{[2]_{p} [2]_{p q}}\\[14pt]
%\noalign{\vskip\alignskip}
\dfrac{-[3]_{p} [2]_{p^2 q^2} \theta}{[2]_{p} [2]_{p/q} [2]_{p^2/q} \eta}&
-g_{25}(-1/p,-1/q)&
\dfrac{-2 [3]_{p} [3]_{q} \theta}{[2]_{p} [2]_{p/q}}\\[14pt]
%\noalign{\vskip\alignskip}
\dfrac{-[2]_{p^2 q^2}}{[2]_{p q} [2]_{p/q} [2]_{p^2/q} [3]_{q} \eta}&
\dfrac{-[2]_{p^2/q^2}}{[2]_{p q} [2]_{p/q} [2]_{p^2 q} [3]_{q} \theta}&
-f_{25}(p,-1/q)
\end{matrix}\right),
\end{multline*}

and where
\begin{align*}
f_{25}(x,y) &= -\dfrac{x^4 y^2+x^2-x^2 y^4+y^2}{x^2 y^4 [2]_{x y} [2]_{x/y}},
\quad\text{and}\displaybreak[0]\\
g_{25}(x,y) &= -\dfrac{x^6 y^4-x^4 y^6-x^4 y^2+x^4+x^4 y^4+x^2 y^2+
          x^2-x^2 y^6+y^2-y^6}{x^4 y^4 [2]_{x} [2]_{x/y} [2]_{x^2 y}}.
\end{align*}

\bibliographystyle{plain}
\bibliography{hecke}
\end{document}

%% file: table.tex
\input nine

\input eight

\catcode `\!=11
\catcode `\@=11

\let\!tacr=\\ % Save meaning of \\ (Needed if TABLE is used with LaTeX

\newdimen\LineThicknessUnit 
\newdimen\StrutUnit            
\newskip \InterColumnSpaceUnit  
\newdimen\ColumnWidthUnit     
\newdimen\KernUnit

\let\!taLTU=\LineThicknessUnit % Used in preamble
\let\!taCWU=\ColumnWidthUnit   % Used in preamble
\let\!taKU =\KernUnit          % Used in preamble

\newtoks\NormalTLTU
\newtoks\NormalTSU
\newtoks\NormalTICSU
\newtoks\NormalTCWU
\newtoks\NormalTKU

\NormalTLTU={1in \divide \LineThicknessUnit by 300 }
\NormalTSU ={\normalbaselineskip
  \divide \StrutUnit by 11 }  % 11 = 8+3 = NormalT Height+Depth Factors
\NormalTICSU={.5em plus 1fil minus .25em}  % .5em = width of a digit
\NormalTCWU ={.5em}
\NormalTKU  ={.5em}

\def\NormalTableUnits{%
  \LineThicknessUnit   =\the\NormalTLTU
  \StrutUnit           =\the\NormalTSU
  \InterColumnSpaceUnit=\the\NormalTICSU
  \ColumnWidthUnit     =\the\NormalTCWU
  \KernUnit            =\the\NormalTKU}
 
\NormalTableUnits

\newcount\LineThicknessFactor    
\newcount\StrutHeightFactor      
\newcount\StrutDepthFactor       
\newcount\InterColumnSpaceFactor 
\newcount\ColumnWidthFactor      
\newcount\KernFactor
\newcount\VspaceFactor

\newcount\TracingKeys % >=1 reports new keys, >=2 reports key usage
\newcount\TracingFormats  % >=1 reports templates for columns
\def\BeginTableParBox#1{%
  \vtop\bgroup 
    \hsize=#1
    \normalbaselines 
    \let~=\!ttTie
    \let\-=\!ttDH
    \the\EveryTableParBox} 
  
\def\EndTableParBox{%
    \MakeStrut{0pt}{\StrutDepthFactor\StrutUnit}
  \egroup} % finishes the \vtop begun by \BeginTableParbox

\newtoks\EveryTableParBox
\EveryTableParBox={%
  \parindent=0pt
  \raggedright
  \rightskip=0pt plus 4em %   Provide more stretch
  \relax}

\newtoks\EveryTable
\newtoks\!taTableSpread

\newskip\LeftTabskip
\newskip\RightTabskip

\newcount\!taCountA
\newcount\!taColumnNumber
\newcount\!taRecursionLevel % (Initially 0)

\newdimen\!taDimenA  % used by \Enlarge
\newdimen\!taDimenB  % used by \Enlarge
\newdimen\!taDimenC  % used by numeric.tex
\newdimen\!taMinimumColumnWidth

\newtoks\!taToksA

\newtoks\!taPreamble
\newtoks\!taDataColumnTemplate
\newtoks\!taRuleColumnTemplate
\newtoks\!taOldRuleColumnTemplate
\newtoks\!taLeftGlue
\newtoks\!taRightGlue

\newskip\!taLastRegularTabskip

\newif\if!taDigit
\newif\if!taBeginFormat
\newif\if!taOnceOnlyTabskip

\def\TaBlE{%
  T\kern-.27em\lower.5ex\hbox{A}\kern-.18em B\kern-.1em
    \lower.5ex\hbox{L}\kern-.075em E}

{\catcode`\|=13 \catcode`\"=13
  \gdef\ActivateBarAndQuote{%
    \ifnum \catcode`\|=13
    \else
      \catcode`\|=13
      \def|{%
        \ifmmode
          \vert
        \else
          \char`\|
        \fi}%
    \fi
    \ifnum \catcode`\"=13
    \else
      \catcode`\"=13
      \def"{\char`\"}%
    \fi}}
 
{\catcode `\|=12 \catcode `\"=12 

}

\def\!thMessage#1{\immediate\write16{#1}\ignorespaces}
 
\let\!thx=\expandafter

\def\!thGobble#1{} 

\def\\{\let\!thSpaceToken= }\\ 

\def\!thHeight{height}
\def\!thDepth{depth}
\def\!thWidth{width}

\def\!thToksEdef#1=#2{%
  \edef\!ttemp{#2}%
  #1\!thx{\!ttemp}%
  \ignorespaces}

\def\!thStoreErrorMsg#1#2{%
  \toks0 =\!thx{\csname #2\endcsname}%
  \edef#1{\the\toks0 }}

\def\!thReadErrorMsg#1{%
  \!thx\!thx\!thx\!thGobble\!thx\string #1}

\def\!thError#1#2{%
  \begingroup
    \newlinechar=`\^^J%
    \edef\!ttemp{#2}%
    \errhelp=\!thx{\!ttemp}%
    \!thMessage{%
      ^^J\!thReadErrorMsg\!thErrorMsgA 
      ^^J\!thReadErrorMsg\!thErrorMsgB}%
    \errmessage{#1}%
  \endgroup}

\!thStoreErrorMsg\!thErrorMsgA{%
  TABLE error; see manual for explanation.}
\!thStoreErrorMsg\!thErrorMsgB{%
  Type \space H <return> \space for immediate help.}

\def\!thGetReplacement#1#2{%
   \begingroup
     \!thMessage{#1}
     \endlinechar=-1
     \global\read16 to#2%
   \endgroup}

\def\!thLoop#1\repeat{%
  \def\!thIterate{%
    #1%
    \!thx \!thIterate
    \fi}%
  \!thIterate 
  \let\!thIterate\relax}

\def\Smash{%
  \relax
  \ifmmode
    \expandafter\mathpalette
    \expandafter\!thDoMathVCS
  \else
    \expandafter\!thDoVCS
  \fi}
                      
\def\!thDoVCS#1{%
  \setbox\z@\hbox{#1}%
  \!thFinishVCS}
                      
\def\!thDoMathVCS#1#2{%
  \setbox\z@\hbox{$\m@th#1{#2}$}%
  \!thFinishVCS}
                      
\def\!thFinishVCS{%
  \vbox to\z@{\vss\box\z@\vss}}

\def\!thSetDimen{%
  \ifnum \!tgCode=1
    \ifx \!tgValue\empty
      \!taDimenA \StrutHeightFactor\StrutUnit
      \advance \!taDimenA \StrutDepthFactor\StrutUnit
      \divide \!taDimenA 2
    \else
      \!taDimenA \!tgValue\StrutUnit
    \fi
  \else
    \!taDimenA \!tgValue
  \fi
  \!taDimenA=\!thSign\!taDimenA\relax
  \ifmmode
    \expandafter\mathpalette
    \expandafter\!thDoMathRaise
  \else
    \expandafter\!thDoSimpleRaise
  \fi}
                      
\def\!thDoSimpleRaise#1{%
  \setbox\z@\hbox{\raise \!taDimenA\hbox{#1}}%
  \!thFinishRaise} % From Plain TeX: \ht0=0pt \dp0=0pt \box0
                      
\def\!thDoMathRaise#1#2{%
  \setbox\z@\hbox{\raise \!taDimenA\hbox{$\m@th#1{#2}$}}%
  \!thFinishRaise}

\def\!thFinishRaise{%
  \ht\z@\z@ 
  \dp\z@\z@
  \box\z@}

\def\!thKernBack{%
  \kern -
  \ifnum \!tgCode=1 
    \ifx \!tgValue\empty 
      \the\KernFactor
    \else
      \!tgValue    % user-specified integer
    \fi
    \KernUnit
  \else 
    \!tgValue      % user-specified dimension
  \fi
  \ignorespaces}%

\def\Vspace{%
  \noalign
  \bgroup
  \!tgGetValue\!thVspace}

\def\!thVspace{%
  \vskip
    \ifnum \!tgCode=1 
      \ifx \!tgValue\empty 
        \the\VspaceFactor
      \else
        \!tgValue    % user-specified integer
      \fi
      \StrutUnit
    \else 
      \!tgValue      % user-specified skip
    \fi
  \egroup} % Ends the \noalign

\def\BeginFormat{%
  \catcode`\|=12 % Inhibit expansion if | immediately follows a <number>
  \catcode`\"=12 %  read by \getvalue.
  \!taPreamble={}% 
  \!taColumnNumber=0
  \skip0 =\InterColumnSpaceUnit
  \multiply\skip0 \InterColumnSpaceFactor
  \divide\skip0 2
  \!taRuleColumnTemplate=\!thx{%
    \!thx\tabskip\the\skip0 }%
  \!taLastRegularTabskip=\skip0 
  \!taOnceOnlyTabskipfalse
  \!taBeginFormattrue % Used to intercept key "]"
  \def\!tfRowOfWidths{}%  Artificial Table Row with horizontal struts
  \ReadFormatKeys}

\def\!tfSetWidth{%
  \ifx \!tfRowOfWidths \empty  % true if no prior "w" keys
    \ifnum \!taColumnNumber>0  % true if "w" key is to right of first "|"
      \begingroup              % RowOfWidths={&\omit || n copies of
         \!taCountA=1          % to the left of this one
         \aftergroup \edef \aftergroup \!tfRowOfWidths \aftergroup {%
           \aftergroup &\aftergroup \omit
           \!thLoop
             \ifnum \!taCountA<\!taColumnNumber
             \advance\!taCountA 1
             \aftergroup \!tfAOAO
           \repeat 
           \aftergroup }%
      \endgroup
    \fi
  \fi      
  \ifx [\!ttemp % \!tgGetValue sets \!ttemp = token after w
    \!thx\!tfSetWidthText
  \else
    \!thx\!tfSetWidthValue
  \fi}

\def\!tfAOAO{%
  &\omit&\omit}

\def\!tfSetWidthText [#1]{% #1 = specified text 
  \def\!tfWidthText{#1}%
  \ReadFormatKeys}

\def\!tfSetWidthValue{%
  \!taMinimumColumnWidth = 
    \ifnum \!tgCode=1 
      \ifx\!tgValue\empty % Use default multiplier if user didn't specify one
        \ColumnWidthFactor
      \else
        \!tgValue 
      \fi
      \ColumnWidthUnit
    \else
      \!tgValue 
    \fi
  \def\!tfWidthText{}%      Override possible prior `w[sample entry]'
  \ReadFormatKeys}

\def\!tfSetTabskip{%
  \ifnum \!tgCode=1
    \skip0 =\InterColumnSpaceUnit
    \multiply\skip0 
      \ifx \!tgValue\empty
        \InterColumnSpaceFactor         % Default integer
      \else
       \!tgValue                        % User-specified integer
      \fi
  \else
    \skip0 =\!tgValue                   % User-specified <skip>
  \fi
  \divide\skip0 by 2
  \ifnum\!taColumnNumber=0 
    \!thToksEdef\!taRuleColumnTemplate={%
      \the\!taRuleColumnTemplate 
      \tabskip \the\skip0 }
  \else
    \!thToksEdef\!taDataColumnTemplate={%
      \the\!taDataColumnTemplate 
      \tabskip \the\skip0 }
  \fi
  \if!taOnceOnlyTabskip
  \else
    \!taLastRegularTabskip=\skip0 % Remember this Tabskip, for possible
  \fi                             % restoration after a subsequent"OnceOnly"
  \ReadFormatKeys}

\def\!tfSetVrule{%
  \!thToksEdef\!taRuleColumnTemplate={%
    \noexpand\hfil
    \noexpand\vrule
    \noexpand\!thWidth
    \ifnum \!tgCode=1
      \ifx \!tgValue\empty
        \the\LineThicknessFactor      % Default integer
      \else
        \!tgValue                     % User-specified integer
      \fi
      \!taLTU                         % \LineThicknessUnit
    \else
      \!tgValue                       % User-specified dimension
    \fi
    ####%
    \noexpand\hfil
    \the\!taRuleColumnTemplate}       % has \tabskips, when column number=0
  \!tfAdjoinPriorColumn}
 
\def\!tfSetAlternateVrule{%
  \afterassignment\!tfSetAlternateA
  \toks0 =}                           % Put template into \toks0

\def\!tfSetAlternateA{%
  \!thToksEdef\!taRuleColumnTemplate={%
    \the\toks0 \the\!taRuleColumnTemplate} % RCT may have \tabskips
  \!tfAdjoinPriorColumn}

\def\!tfAdjoinPriorColumn{%
  \ifnum \!taColumnNumber=0
    \!taPreamble=\!taRuleColumnTemplate % New \tabskip may have been added
    \ifnum \TracingFormats>0             
      \!tfShowRuleTemplate
    \fi
  \else
    \ifx\!tfRowOfWidths\empty  % no "w" keys specified yet, not even this col
    \else
      \!tfUpdateRowOfWidths
    \fi
    \!thToksEdef\!taDataColumnTemplate={%
      \the \!taLeftGlue
      \the \!taDataColumnTemplate
      \the \!taRightGlue}
    \ifnum \TracingFormats>0
      \!tfShowTemplates
    \fi
    \!thToksEdef\!taPreamble={%
      \the\!taPreamble
      &
      \the\!taDataColumnTemplate
      &
      \the\!taRuleColumnTemplate}
  \fi
  \advance \!taColumnNumber 1
  \if!taOnceOnlyTabskip              
    \!thToksEdef\!taDataColumnTemplate={%
       ####\tabskip \the\!taLastRegularTabskip}
  \else
    \!taDataColumnTemplate{##}%
  \fi
  \!taRuleColumnTemplate{}% # is inserted by \SetVrule, or \SetAlternateVrule
  \!taLeftGlue{\hfil}%         % Default positioning is "center"
  \!taRightGlue{\hfil}%         
  \!taMinimumColumnWidth=0pt
  \def\!tfWidthText{}%
  \!taOnceOnlyTabskipfalse    % Set true by key "o"
  \ReadFormatKeys}

\def\!tfUpdateRowOfWidths{%
  \ifx \!tfWidthText\empty
  \else % set specified text according to current template & find width
    \!tfComputeMinColWidth
  \fi
  \edef\!tfRowOfWidths{%
    \!tfRowOfWidths
    &%
    \omit                                  % Data Column
    \ifdim \!taMinimumColumnWidth>0pt
      \hskip \the\!taMinimumColumnWidth
    \fi
    &
    \omit}}                                % Rule Column

\def\!tfComputeMinColWidth{%
  \setbox0 =\vbox{%
    \ialign{% Plain's initialized \halign; \tabskip=0pt \everycr={}
       \span\the\!taDataColumnTemplate\cr
       \!tfWidthText\cr}}%
  \!taMinimumColumnWidth=\wd0 }

\def\!tfShowRuleTemplate{%
  \!thMessage{}
  \!thMessage{TABLE FORMAT}
  \!thMessage{Column: Template}
  \!thMessage{%
    \space *c: ##\tabskip \the\LeftTabskip}
  \!taOldRuleColumnTemplate=\!taRuleColumnTemplate}

\def\!tfShowTemplates{%
  \!thMessage{%
    \space \space r: \the\!taOldRuleColumnTemplate}
  \!taOldRuleColumnTemplate=\!taRuleColumnTemplate
  \!thMessage{%
    \ifnum \!taColumnNumber<10
      \space
    \fi
    \the\!taColumnNumber c: \the\!taDataColumnTemplate}
  \ifdim\!taMinimumColumnWidth>0pt
    \!thMessage{%
      \space \space w: \the\!taMinimumColumnWidth}
  \fi}

\def\!tfFinishFormat{%
  \ifnum \TracingFormats>0
    \!thMessage{%
      \space \space r: \the\!taOldRuleColumnTemplate
        \tabskip \the\RightTabskip}%
    \!thMessage{%
      \space *c: ##\tabskip 0pt}
  \fi
  \ifnum \!taColumnNumber<2
    \!thError{%
      \ifnum \!taColumnNumber=0
        No
      \else
        Only 1
      \fi
      "|"}%
      {\!thReadErrorMsg\!tfTooFewBarsA
       ^^J\!thReadErrorMsg\!tfTooFewBarsB
       ^^J\!thReadErrorMsg\!tkFixIt}%
  \fi
  \!thToksEdef\!taPreamble={%
    ####\tabskip\LeftTabskip 
    &
    \the\!taPreamble \tabskip\RightTabskip
    &
    ####\tabskip 0pt \cr}
  \ifnum \TracingFormats>1
    \!thMessage{Preamble=\the\!taPreamble}
  \fi
  \ifnum \TracingFormats>2
    \!thMessage{Row Of Widths="\!tfRowOfWidths"}
  \fi
  \!taBeginFormatfalse % Intercepts "|", tabskips, and "."
  \catcode`\|=13
  \catcode`\"=13
  \!ttDoHalign}

\!thStoreErrorMsg\!tfTooFewBarsA{%
  There must be at least 2 "|"'s (and/or "\string \|"'s)}
\!thStoreErrorMsg\!tfTooFewBarsB{%
  between \string\BeginFormat\space and \string\EndFormat\space (or ".").}

\def\ReFormat[{%
  \omit
  \!taDataColumnTemplate{##}%
  \!taLeftGlue{}% 
  \!taRightGlue{}% 
  \catcode`\|=12  % Inhibit expansion if | immediately follows a <number>
  \catcode`\"=12  % read by \getvalue. Actually, '|' and '"' shouldn't
  \ReadFormatKeys}% appear in a \ReFormat cmd; this is here as a safeguard.

\def\!tfEndReFormat{%
  \ifnum \TracingFormats>0
    \!thMessage{ReF: 
       \the\!taLeftGlue
       \hbox{\the\!taDataColumnTemplate}%   White lie
       \the\!taRightGlue}
  \fi
  \catcode`\|=13
  \catcode`\"=13
  \!tfReFormat}

\def\!tfReFormat#1{%
  \the \!taLeftGlue
  \vbox{%
    \ialign{%
      \span\the\!taDataColumnTemplate\cr
       #1\cr}}%
  \the \!taRightGlue}

\def\!tgGetValue#1{%
  \def\!tgReturn{#1}%          Set return
  \futurelet\!ttemp\!tgCheckForParen}%  Now \!ttemp is the token 
\def\!tgCheckForParen{%
  \ifx\!ttemp (%
    \!thx \!tgDoParen
  \else
    \!thx \!tgCheckForSpace
  \fi}

\def\!tgDoParen(#1){%
  \def\!tgCode{2}%
  \def\!tgValue{#1}%     NOTE #1 MUST BE A LEGITIMATE VALUE
  \!tgReturn}

\def\!tgCheckForSpace{%
  \def\!tgCode{1}%
  \def\!tgValue{}%           Initialize value to <null>
  \ifx\!ttemp\!thSpaceToken
    \!thx \!tgReturn        % <blank space> means no value was specified
  \else
    \!thx \!tgCheckForDigit         
  \fi}

\def\!tgCheckForDigit{%
  \!taDigitfalse
  \ifx 0\!ttemp
    \!taDigittrue
  \else
    \ifx 1\!ttemp
      \!taDigittrue
    \else
      \ifx 2\!ttemp
        \!taDigittrue
      \else
        \ifx 3\!ttemp
          \!taDigittrue
        \else
          \ifx 4\!ttemp
            \!taDigittrue
          \else
            \ifx 5\!ttemp
              \!taDigittrue
            \else
              \ifx 6\!ttemp
                \!taDigittrue
              \else
                \ifx 7\!ttemp
                  \!taDigittrue
                \else
                  \ifx 8\!ttemp
                    \!taDigittrue
                  \else
                    \ifx 9\!ttemp
                      \!taDigittrue
                    \fi
                  \fi
                \fi
              \fi
            \fi
          \fi
        \fi
      \fi
    \fi
  \fi
  \if!taDigit
    \!thx \!tgGetNumber
  \else
    \!thx \!tgReturn 
  \fi}

\def\!tgGetNumber{%
  \afterassignment\!tgGetNumberA
  \!taCountA=}
\def\!tgGetNumberA{%
  \edef\!tgValue{\the\!taCountA}%
  \!tgReturn}

\def\!tgSetUpParBox{%
  \edef\!ttemp{%
    \noexpand \ReadFormatKeys
    b{\noexpand \BeginTableParBox{%
      \ifnum \!tgCode=1 
        \ifx \!tgValue\empty 
          \the\ColumnWidthFactor
        \else
          \!tgValue    % user-specified integer
        \fi
        \!taCWU        % \ColumnWidthUnit 
      \else 
        \!tgValue      % user-specified dimension
      \fi}}}%
  \!ttemp
  a{\EndTableParBox}}

\def\!tgInsertKern{%
  \edef\!ttemp{%
    \kern
    \ifnum \!tgCode=1 
      \ifx \!tgValue\empty 
        \the\KernFactor
      \else
        \!tgValue    % user-specified integer
      \fi
      \!taKU         % \KernUnit
    \else 
      \!tgValue      % user-specified dimension
    \fi}%
  \edef\!ttemp{%
    \noexpand\ReadFormatKeys
    \ifh@            % true if kern goes to left of "#"
      b{\!ttemp}
    \fi
    \ifv@            % true if kern goes to right of "#"
      a{\!ttemp}
    \fi}%
  \!ttemp}

\def\NewFormatKey#1{%
  \!thx\def\!thx\!ttempa\!thx{\string #1}%
  \!thx\def\!thx\!ttempb\!thx{\csname !tk<\!ttempa>\endcsname}%
  \ifnum \TracingKeys>0
    \!tkReportNewKey
  \fi
  \!thx\ifx \!ttempb \relax
    \!thx\!tkDefineKey
  \else 
    \!thx\!tkRejectKey
  \fi}

\def\!tkReportNewKey{%
  \!taToksA\!thx{\!ttempa}%  
  \!thMessage{NEW KEY: "\the\!taToksA"}}

\def\!tkDefineKey{%
  \!thx\def\!ttempb}%

\def\!tkRejectKey{%
    \!taToksA\!thx{\!ttempa}% 
    \!thError{Key letter "\the\!taToksA" already used}
      {\!thReadErrorMsg\!tkFixIt}
    \def\!tkGarbage}%

\!thStoreErrorMsg\!tkFixIt{%
  You'd better type \space 'E' \space and fix your file.}

\def\ReadFormatKeys#1{%
  \!thx\def\!thx\!ttempa\!thx{\string #1}%
  \!thx\def\!thx\!ttempb\!thx{\csname !tk<\!ttempa>\endcsname}%
  \ifnum \TracingKeys>1
    \!tkReportKey
  \fi
  \!thx\ifx \!ttempb\relax 
    \!thx\!tkReplaceKey
  \else
    \!thx\!ttempb
  \fi}

\def\!tkReportKey{%
  \!taToksA\!thx{\!ttempa}% 
  \!thMessage{KEY: "\the\!taToksA"}}

\def\!tkReplaceKey{%
  \!taToksA\!thx{\!ttempa}%
  \!thError {Undefined format key "\the\!taToksA"}
    {\!thReadErrorMsg\!tkUndefined ^^J\!thReadErrorMsg\!tkBadKey}
  \!tkReplaceKeyA}

\def\!tkReplaceKeyA{%
  \!thGetReplacement{\!thReadErrorMsg\!tkReplace}\!tkReplacement
  \!thx\ReadFormatKeys\!tkReplacement}

\!thStoreErrorMsg\!tkUndefined{%
  The format key in " "'s on the next to top line is undefined.}
\!thStoreErrorMsg\!tkBadKey{%
  Type \space E \space to quit now, or
  \space<CR> \space and respond to next prompt.}
\!thStoreErrorMsg\!tkReplace{%
  Type \space<replacement key><CR> \space,
   or simply \space<CR> \space to skip offending key:}

\NewFormatKey b#1{%
  \!thx\!tkJoin\!thx{\the\!taDataColumnTemplate}{#1}%
  \ReadFormatKeys}

\def\!tkJoin#1#2{%
  \!taDataColumnTemplate{#2#1}}%

\NewFormatKey a#1{%
  \!taDataColumnTemplate\!thx{\the\!taDataColumnTemplate #1}%
  \ReadFormatKeys}

\NewFormatKey \{{%
  \!taDataColumnTemplate=\!thx{\!thx{\the\!taDataColumnTemplate}}%
  \ReadFormatKeys}

\NewFormatKey *#1#2{%
  \!taCountA=#1\relax
  \!taToksA={}%
  \!thLoop 
    \ifnum \!taCountA > 0
    \!taToksA\!thx{\the\!taToksA #2}%
    \advance\!taCountA -1
  \repeat 
  \!thx\ReadFormatKeys\the\!taToksA}

\NewFormatKey \LeftGlue#1{%
  \!taLeftGlue{#1}%
  \ReadFormatKeys}

\NewFormatKey \RightGlue#1{%
  \!taRightGlue{#1}%
  \ReadFormatKeys}

\NewFormatKey c{%
  \ReadFormatKeys 
  \LeftGlue\hfil
  \RightGlue\hfil}

\NewFormatKey l{%
  \ReadFormatKeys 
  \LeftGlue{}   % In case more than one positioning key is specified.
  \RightGlue\hfil}

\NewFormatKey r{%
  \ReadFormatKeys 
  \LeftGlue\hfil
  \RightGlue{}}

\NewFormatKey k{%
  \h@true
  \v@true
  \!tgGetValue{\!tgInsertKern}}

\NewFormatKey i{%
  \h@true
  \v@false
  \!tgGetValue{\!tgInsertKern}}
  
\NewFormatKey j{%
  \h@false
  \v@true
  \!tgGetValue{\!tgInsertKern}}

\NewFormatKey n{%
  \def\!tnStyle{}%
   \futurelet\!tnext\!tnTestForBracket}

\NewFormatKey N{%
  \def\!tnStyle{$}%
   \futurelet\!tnext\!tnTestForBracket}

\NewFormatKey m{%
  \ReadFormatKeys b$ a$}

\NewFormatKey M{%
  \ReadFormatKeys \{ b{$\displaystyle} a$}

\NewFormatKey \m{%
  \ReadFormatKeys l b{{}} m}

\NewFormatKey \M{%
  \ReadFormatKeys l b{{}} M}

\NewFormatKey f#1{%
  \ReadFormatKeys b{#1}}

\NewFormatKey B{%
  \ReadFormatKeys f\bf}

\NewFormatKey I{%
  \ReadFormatKeys f\it}

\NewFormatKey S{%
  \ReadFormatKeys f\sl}

\NewFormatKey R{%
  \ReadFormatKeys f\rm}

\NewFormatKey T{%
  \ReadFormatKeys f\tt}

\NewFormatKey p{%
  \!tgGetValue{\!tgSetUpParBox}}

\NewFormatKey w{%
  \!tkTestForBeginFormat w{\!tgGetValue{\!tfSetWidth}}}

\NewFormatKey s{%
  \!taOnceOnlyTabskipfalse    % in case same column has a prior "o" key
  \!tkTestForBeginFormat t{\!tgGetValue{\!tfSetTabskip}}}

\NewFormatKey o{%
  \!taOnceOnlyTabskiptrue
  \!tkTestForBeginFormat o{\!tgGetValue{\!tfSetTabskip}}}

\NewFormatKey |{%
  \!tkTestForBeginFormat |{\!tgGetValue{\!tfSetVrule}}}

\NewFormatKey \|{%
  \!tkTestForBeginFormat \|{\!tfSetAlternateVrule}}

\NewFormatKey .{%
  \!tkTestForBeginFormat.{\!tfFinishFormat}} 

\NewFormatKey \EndFormat{%
  \!tkTestForBeginFormat\EndFormat{\!tfFinishFormat}} 

\NewFormatKey ]{%
  \!tkTestForReFormat ] \!tfEndReFormat}

\def\!tkTestForBeginFormat#1#2{% 
  \if!taBeginFormat  
    \def\!ttemp{#2}%
    \!thx \!ttemp    
  \else
    \toks0={#1}%  
    \toks2=\!thx{\string\ReFormat}%
    \!thx \!tkImproperUse
  \fi}   

\def\!tkTestForReFormat#1#2{% 
  \if!taBeginFormat  
    \toks0={#1}%  
    \toks2=\!thx{\string\BeginFormat}%
    \!thx \!tkImproperUse
  \else
    \def\!ttemp{#2}%
    \!thx \!ttemp    
  \fi}   

\def\!tkImproperUse{% 
  \!thError{\!thReadErrorMsg\!tkBadUseA "\the\toks0 "}%
    {\!thReadErrorMsg\!tkBadUseB \the\toks2 \space command.
    ^^J\!thReadErrorMsg\!tkBadKey}%
  \!tkReplaceKeyA}
 
\!thStoreErrorMsg\!tkBadUseA{Improper use of key }  
\!thStoreErrorMsg\!tkBadUseB{% 
  The key mentioned above can't be used in a }

\def\!tnTestForBracket{%
  \ifx [\!tnext
    \!thx\!tnGetArgument
  \else
    \!thx\!tnGetCode
  \fi}

\def\!tnGetCode#1 {%  NOTE THE BLANK
  \!tnConvertCode #1..!}

\def\!tnConvertCode #1.#2.#3!{%
  \begingroup
    \aftergroup\edef \aftergroup\!ttemp \aftergroup{%
      \aftergroup[% 
      \!taCountA #1
      \!thLoop
        \ifnum \!taCountA>0
        \advance\!taCountA -1
        \aftergroup0
      \repeat
      \def\!ttemp{#3}%
      \ifx\!ttemp \empty
      \else
        \aftergroup.
        \!taCountA #2
        \!thLoop 
          \ifnum \!taCountA>0
          \advance\!taCountA -1
          \aftergroup0
        \repeat
      \fi 
      \aftergroup]\aftergroup}%
    \endgroup\relax
    \!thx\!tnGetArgument\!ttemp}
  
\def\!tnGetArgument[#1]{%
  \!tnMakeNumericTemplate\!tnStyle#1..!}

\def\!tnMakeNumericTemplate#1#2.#3.#4!{%  #1=<empty> or $
  \def\!ttemp{#4}%
  \ifx\!ttemp\empty
    \!taDimenC=0pt
  \else
    \setbox0=\hbox{\m@th #1.#3#1}%
    \!taDimenC=\wd0
  \fi
  \setbox0 =\hbox{\m@th #1#2#1}%
  \!thToksEdef\!taDataColumnTemplate={%
    \noexpand\!tnSetNumericItem
    {\the\wd0 }%
    {\the\!taDimenC}%
    {#1}%
    \the\!taDataColumnTemplate}  % Might have tabskip glue in here
  \ReadFormatKeys}

\def\!tnSetNumericItem #1#2#3#4 {%             NOTE THE BLANK
  \!tnSetNumericItemA {#1}{#2}{#3}#4..!}

\def\!tnSetNumericItemA #1#2#3#4.#5.#6!{%
  \def\!ttemp{#6}%
  \hbox to #1{\hss \m@th #3#4#3}%
  \hbox to #2{%
    \ifx\!ttemp\empty
    \else
       \m@th #3.#5#3%
    \fi
    \hss}}

\def\MakeStrut#1#2{%
  \vrule width0pt height #1 depth #2}

\def\StandardTableStrut{%
  \MakeStrut{\StrutHeightFactor\StrutUnit}
    {\StrutDepthFactor\StrutUnit}}

\def\AugmentedTableStrut#1#2{%
  \dimen@=\StrutHeightFactor\StrutUnit
  \advance\dimen@ #1\StrutUnit
  \dimen@ii=\StrutDepthFactor\StrutUnit
  \advance\dimen@ii #2\StrutUnit
  \MakeStrut{\dimen@}{\dimen@ii}}

\def\Enlarge#1#2{%  3rd argument is picked up later
  \!taDimenA=#1\relax
  \!taDimenB=#2\relax
  \let\!TsSpaceFactor=\empty
  \ifmmode
    \!thx \mathpalette
    \!thx \!TsEnlargeMath
  \else
    \!thx \!TsEnlargeOther
  \fi}

\def\!TsEnlargeOther#1{%
  \ifhmode
    \setbox\z@=\hbox{#1%
      \xdef\!TsSpaceFactor{\spacefactor=\the\spacefactor}}%
  \else
    \setbox\z@=\hbox{#1}%
  \fi
  \!TsFinishEnlarge}
    
\def\!TsEnlargeMath#1#2{%
  \setbox\z@=\hbox{$\m@th#1{#2}$}%
  \!TsFinishEnlarge}

\def\!TsFinishEnlarge{%
  \dimen@=\ht\z@
  \advance \dimen@ \!taDimenA
  \ht\z@=\dimen@
  \dimen@=\dp\z@
  \advance \dimen@ \!taDimenB
  \dp\z@=\dimen@
  \box\z@ \!TsSpaceFactor{}}

\def\OpenUp#1#2{%
  \advance \StrutHeightFactor #1\relax
  \advance \StrutDepthFactor #2\relax}

\def\BeginTable{%
  \futurelet\!tnext\!ttBeginTable}

\def\!ttBeginTable{%
  \ifx [\!tnext
    \def\!tnext{\!ttBeginTableA}%
  \else 
    \def\!tnext{\!ttBeginTableA[c]}%
  \fi
  \!tnext}

\def\!ttBeginTableA[#1]{%
  \if #1u%                  % "unboxed" table
    \ifmmode                 
      \def\!ttEndTable{%    % user had better be in display math mode
        \relax}%            %   and have only one table at the outer level
    \else                   % user had better be in vertical mode
      \bgroup
      \def\!ttEndTable{%
        \egroup}%
    \fi
  \else
    \hbox\bgroup $
    \def\!ttEndTable{%
      \egroup %   for the \vtop, \vbox, or \vcenter, yet to come
      $%          for math mode
      \egroup}%   for the \hbox
    \if #1t%
      \vtop
    \else
      \if #1b%
        \vbox
      \else
        \vcenter % math mode was essential for this
      \fi
    \fi
    \bgroup      % for the \vtop, \vbox, or \vcenter
  \fi
  \advance\!taRecursionLevel 1 % RecursionLevel governs initialization
  \let\!ttRightGlue=\relax  % This may be changed by \JustCenter, etc
  \everycr={}
  \ifnum \!taRecursionLevel=1
    \!ttInitializeTable
  \fi}

\bgroup
  \catcode`\|=13
  \catcode`\"=13
  \catcode`\~=13
  \gdef\!ttInitializeTable{%
    \let\!ttTie=~ %                         Meanings of ~ and \- are
    \let\!ttDH=\- %                           restored by \BeginTableParBox
    \catcode`\|=\active
    \catcode`\"=\active
    \catcode`\~=\active
    \def |{\unskip\!ttRightGlue&&}%         Use rule-column template
    \def\|{\unskip\!ttRightGlue&\omit\!ttAlternateVrule}% 
    \def"{\unskip\!ttRightGlue&\omit&}%     Omit rule-column template
    \def~{\kern .5em}%                      ~ now has the width of a digit
    \def\\{\!ttEndOfRow}%
    \def\-{\!ttShortHrule}%
    \def\={\!ttLongHrule}%
    \def\_{\!ttFullHrule}%
    \def\Left##1{##1\hfill\null}%           \null prevents \unskip from
    \def\Center##1{\hfill ##1\hfill\null}%    killing the \hfill
    \def\Right##1{\hfill##1}%
    \the\EveryTable}
\egroup

\let\!ttRightGlue=\relax  % This may be changed, in a group, by 
\def\!ttDoHalign{%
  \baselineskip=0pt \lineskiplimit=0pt \lineskip=0pt %
  \tabskip=0pt
  \halign \the\!taTableSpread \bgroup
   \span\the\!taPreamble
   \ifx \!tfRowOfWidths \empty
   \else 
     \!tfRowOfWidths \cr % 
   \fi}

\def\EndTable{%
  \egroup % finishes the \halign
  \!ttEndTable}%    closes off the table envirnoment set up by \BeginTable

\def\!ttEndOfRow{%
  \futurelet\!tnext\!ttTestForBlank}

\def\!ttTestForBlank{%
  \ifx \!tnext\!thSpaceToken  % the "usual" case
    \!thx\!ttDoStandard
  \else
    \!thx\!ttTestForZero
  \fi}
  
\def\!ttTestForZero{%
  \ifx 0\!tnext
    \!thx \!ttDoZero
  \else
    \!thx \!ttTestForPlus
  \fi}

\def\!ttTestForPlus{%
  \ifx +\!tnext
    \!thx \!ttDoPlus
  \else
    \!thx \!ttDoStandard
  \fi}

\def\!ttDoZero#1{% #1 eats the 0
  \cr} 

\def\!ttDoPlus#1#2#3{% #1 eats the +
  \AugmentedTableStrut{#2}{#3}%
  \cr} 

\def\!ttDoStandard{% 
  \StandardTableStrut
  \cr}

\def\!ttAlternateVrule{%
  \!tgGetValue{\!ttAVTestForCode}}  % AV == Alternate Vrule

\def\!ttAVTestForCode{%
  \ifnum \!tgCode=2              % (...) follows "\|"
    \!thx\!ttInsertVrule         % \InsertVrule ends with "&"
  \else
    \!thx\!ttAVTestForEmpty
  \fi}

\def\!ttAVTestForEmpty{%
  \ifx \!tgValue\empty           % non-digit after "\|"
    \!thx\!ttAVTestForBlank
  \else
    \!thx\!ttInsertVrule         % integer after "\|"
  \fi}

\def\!ttAVTestForBlank{%
  \ifx \!ttemp\!thSpaceToken     % blank after "\|"
    \!thx\!ttInsertVrule
  \else
    \!thx\!ttAVTestForStar 
  \fi}

\def\!ttAVTestForStar{%
  \ifx *\!ttemp                  % "*" after "\|"
    \!thx\!ttInsertDefaultPR     % PR == pseudo-rule
  \else
    \!thx\!ttGetPseudoVrule       % "Anything else" after "\|"
  \fi}

\def\!ttInsertVrule{%
  \hfil 
  \vrule \!thWidth
    \ifnum \!tgCode=1
      \ifx \!tgValue\empty 
        \LineThicknessFactor
      \else
        \!tgValue
      \fi
      \LineThicknessUnit
    \else
      \!tgValue
    \fi
  \hfil
  &}

\def\!ttInsertDefaultPR*{%
  \PseudoVrule    % User-specified default pseudo-rule
  &}

\def\!ttGetPseudoVrule#1{%
  \toks0={#1}%
  #1&}

\def\PseudoVrule{}

\def\!ttuse#1{%
  \ifnum #1>\@ne 
    \omit 
    \mscount=#1 %        \mscount is in Plain
    \advance\mscount by \m@ne
    \advance\mscount by \mscount
    \!thLoop 
      \ifnum\mscount>\@ne 
      \sp@n %            from Plain (\span\omit \advance\mscount\m@ne)
    \repeat 
    \span 
  \fi}

\def\!ttUse#1[{%
  \!ttuse{#1}%
  \ReFormat[}

\def\!ttFullHrule{%
  \noalign
  \bgroup
  \!tgGetValue{\!ttFullHruleA}}

\def\!ttFullHruleA{%
  \!ttGetHalfRuleThickness % Sets \dimen0 to half of specified thickness
  \hrule \!thHeight \dimen0 \!thDepth \dimen0
  \penalty0 % so can break an ``unboxed'' table after a horizontal rule.
  \egroup} % ends the \noalign

\def\!ttShortHrule{%
  \omit
  \!tgGetValue{\!ttShortHruleA}}

\def\!ttShortHruleA{%
  \!ttGetHalfRuleThickness % Sets \dimen0 to half of specified thickness
  \leaders \hrule \!thHeight \dimen0 \!thDepth \dimen0 \hfill
  \null    % prevents an \unskip from annihilating the \leaders
  \ignorespaces} 

\def\!ttLongHrule{%
  \omit\span\omit\span \!ttShortHrule}

\def\!ttGetHalfRuleThickness{%
  \dimen0 =
    \ifnum \!tgCode=1
      \ifx \!tgValue\empty
        \LineThicknessFactor
      \else
        \!tgValue    % user-specified integer
      \fi
      \LineThicknessUnit
    \else
      \!tgValue      % user-specified dimension
    \fi
  \divide\dimen0 2 }

\def\WidenTableBy#1{%
  \ifdim #1=0pt
    \!taTableSpread={}%
  \else
    \!taTableSpread={spread #1}%
  \fi}

\def\JustLeft{%
  \omit \let\!ttRightGlue=\hfill}
\def\JustCenter{%
  \omit \hfill\null \let\!ttRightGlue=\hfill}

\let\\=\!tacr
\catcode`\!=12
\catcode`\@=12

%% file: nine.tex
\font\ninerm=cmr9
\font\sevenrm=cmr7
\font\sixrm=cmr6
\font\fiverm=cmr5
\font\ninei=cmmi9
\font\sixi=cmmi6
\font\fivei=cmmi6
\font\ninesy=cmsy9
\font\sixsy=cmsy6
\font\fivesy=cmsy5
\font\tenex=cmex10
\font\nineit=cmti9
\font\ninesl=cmsl9
\font\ninett=cmtt9
\font\ninebf=cmbx9
\font\sixbf=cmbx6
\font\fivebf=cmbx5

\def\ninepoint{\def\rm{\fam0\ninerm}% switch to 9-point type
  \textfont0=\ninerm \scriptfont0=\sixrm \scriptscriptfont0=\fiverm
  \textfont1=\ninei \scriptfont1=\sixi \scriptscriptfont0=\fivei
  \textfont2=\ninesy \scriptfont2=\sixsy \scriptscriptfont2=\fivesy
  \textfont3=\tenex \scriptfont3=\tenex \scriptscriptfont3=\tenex
  \textfont\itfam=\nineit  \def\it{\fam\itfam\nineit}%
  \textfont\slfam=\ninesl  \def\sl{\fam\slfam\ninesl}%
  \textfont\ttfam=\ninett  \def\tt{\fam\ttfam\ninett}%
  \textfont\bffam=\ninebf  \scriptfont\bffam=\sixbf
   \scriptscriptfont\bffam=\fivebf  \def\bf{\fam\bffam\ninebf}%
  \normalbaselineskip=11pt
  \setbox\strutbox=\hbox{\vrule height8pt depth3pt width0pt}%
  \let\sc=\sevenrm  \normalbaselines\rm}

%% file: eight.tex
\font\eightrm=cmr8
\font\sevenrm=cmr7
\font\sixrm=cmr6
\font\fiverm=cmr5
\font\eighti=cmmi8
\font\sixi=cmmi6
\font\fivei=cmmi6
\font\eightsy=cmsy8
\font\sixsy=cmsy6
\font\fivesy=cmsy5
\font\tenex=cmex10
\font\eightit=cmti8
\font\eightsl=cmsl8
\font\eighttt=cmtt8
\font\eightbf=cmbx8
\font\sixbf=cmbx6
\font\fivebf=cmbx5

\def\eightpoint{\def\rm{\fam0\eightrm}% switch to 8-point type
  \textfont0=\eightrm \scriptfont0=\sixrm \scriptscriptfont0=\fiverm
  \textfont1=\eighti \scriptfont1=\sixi \scriptscriptfont0=\fivei
  \textfont2=\eightsy \scriptfont2=\sixsy \scriptscriptfont2=\fivesy
  \textfont3=\tenex \scriptfont3=\tenex \scriptscriptfont3=\tenex
  \textfont\itfam=\eightit  \def\it{\fam\itfam\eightit}%
  \textfont\slfam=\eightsl  \def\sl{\fam\slfam\eightsl}%
  \textfont\ttfam=\eighttt  \def\tt{\fam\ttfam\eighttt}%
  \textfont\bffam=\eightbf  \scriptfont\bffam=\sixbf
   \scriptscriptfont\bffam=\fivebf  \def\bf{\fam\bffam\eightbf}%
  \normalbaselineskip=9pt
  \setbox\strutbox=\hbox{\vrule height7pt depth2pt width0pt}%
  \let\sc=\sixrm  \normalbaselines\rm}

%% file: 1997-027.bbl
\begin{thebibliography}{10}

\bibitem{bourbaki:1968}
N.~Bourbaki.
\newblock {\em Groupes et Alg\`ebres de Lie, Chap. {IV},{V},{VI}}.
\newblock Hermann, Paris, 1968.

\bibitem{carter:1985}
R.~Carter.
\newblock {\em Finite groups of Lie type: Conjugacy classes and complex
  characters}.
\newblock Wiley, New York, 1985.

\bibitem{char-etal:1991}
Bruce~W. Char, Keith~O. Geddes, Gaston~H. Gonnet, Benton~L. Leong, Michael~B.
  Monagan, and Stephen~M. Watt.
\newblock {\em Maple {V} Language Reference Manual}.
\newblock Springer-Verlag, New York, Berlin, Heidelberg, 1991.

\bibitem{curtis:1981}
C.W. Curtis and I.~Reiner.
\newblock {\em Methods of representation theory, Vol. I. With applications to
  finite groups and orders}.
\newblock Wiley, New York, 1981.

\bibitem{geck:1994}
M.~Geck.
\newblock On the character values of {I}wahori-{H}ecke algebras of exceptional
  type.
\newblock {\em Proc. London Math. Soc. (3)}, 68:51--76, 1994.

\bibitem{hoefsmit:1974}
P.N. Hoefsmit.
\newblock {\em Representations of Hecke algebras of finite groups with BN-pairs
  of classical type}.
\newblock Ph.D. Thesis, University of British Columbia, 1974.

\bibitem{kilmoyer:1969}
R.~Kilmoyer.
\newblock {\em Some irreducible representations of a finite group with
  BN-pair}.
\newblock Ph.D. Thesis, Massachusetts Institute of Technology, 1969.

\bibitem{ram:1997}
A.~Ram.
\newblock Seminormal representations of {W}eyl groups and {I}wahori-{H}ecke
  algebras.
\newblock {\em Proc. London Math. Soc. (3)}, 75:999--999, 1997.

\bibitem{wenzl:1988}
H.~Wenzl.
\newblock Hecke algebras of type {$A_n$} and subfactors.
\newblock {\em Invent. Math.}, 92:349--383, 1988.

\bibitem{young:1932}
A.~Young.
\newblock Quantitative substitutional analysis (sixth paper).
\newblock {\em Proc. London Math. Soc. (2)}, 34:196--230, 1932.

\end{thebibliography}
